\documentclass[preprint,authoryear]{elsarticle}
\usepackage[T1]{fontenc}
\usepackage[utf8]{inputenc}
\usepackage{tabularx}
\usepackage{supertabular}
\usepackage{booktabs}
\RequirePackage[authoryear]{natbib}
\usepackage{lmodern}
\usepackage[onehalfspacing]{setspace}
\usepackage{nicefrac}
\usepackage{amsmath}
\usepackage{a4wide}
\usepackage{amssymb}
\usepackage{bbm}
\usepackage{changes}
\usepackage[flushleft]{threeparttable}
\definecolor{dblue}{rgb}{0.21,0.21,0.55}

\usepackage{multirow}
\usepackage{lscape}
\usepackage[english]{babel}
\usepackage{graphicx}
\usepackage{arydshln}
\usepackage{rotating}
\usepackage[flushmargin,hang]{footmisc}
\usepackage{amsmath}
\usepackage{lscape}

\renewcommand{\P}{\mathbb{P}}

\newcommand{\E}{\mathbb{E}}
\newcommand{\N}{\mathbb{N}}
\newcommand{\R}{\mathbb{R}}
\newcommand{\1}{\mathbbm{1}}
\newcommand{\KLEINO}{{\scriptstyle{\mathcal{O}}}}
\DeclareMathAccent{\verywidehat}{\mathord}{largesymbols}{'144}
\newcommand{\var}{\mathbb{V}\hspace*{-0.05cm}\textnormal{a\hspace*{0.02cm}r}}

\newcommand{\cov}{\mathbb{C}\textnormal{o\hspace*{0.02cm}v}}

\allowdisplaybreaks[3]
\newtheorem{remark}{Remark}

\newtheorem{assump}{Assumption}
\newtheorem{prop}{Proposition}[section]

\newtheorem{example}{Example}
\newtheorem{cor}[prop]{Corollary}

\begin{document}
\renewcommand*{\thefootnote}{\fnsymbol{footnote}}

\title{Estimation of the discontinuous leverage effect: Evidence from the NASDAQ order book}
\author[1]{Markus Bibinger}
\author[2]{Christopher Neely}
\author[3]{Lars Winkelmann}
\address[1]{Faculty of Mathematics and Computer Science, Philipps-Universit\"at Marburg} 
\address[2]{Research Department, Federal Reserve Bank of St.\;Louis\footnote[2]{The views expressed are those of the individual authors and do not necessarily reflect official positions of the Federal Reserve Bank of St.\ Louis, the Federal Reserve System, or the Board of Governors.}}
\address[3]{Department of Economics, Freie Universit\"at Berlin}
\normalsize
\begin{frontmatter}
\onehalfspacing
\vspace{-.7cm} 

\begin{abstract}
{{\normalsize \noindent
An extensive empirical literature documents a generally negative correlation, named the ``leverage effect,'' between asset returns and changes of volatility. It is more challenging to establish such a return-volatility relationship for jumps in high-frequency data.
We propose new nonparametric methods to assess and test for a discontinuous leverage effect ---  i.e.\ a relation between contemporaneous jumps in prices and volatility. The methods are robust to market microstructure noise and build on a newly developed price-jump localization and estimation procedure. Our empirical investigation of six years of transaction data from 320 NASDAQ firms displays no unconditional negative correlation between price and volatility cojumps. We show, however, that there is a strong relation between price-volatility cojumps if one conditions on the sign of price jumps and whether the price jumps are market-wide or idiosyncratic. Firms' volatility levels strongly explain the cross-section of discontinuous leverage while debt-to-equity ratios have no significant explanatory power.}}


\begin{keyword}
High-frequency data \sep market microstructure \sep news impact\sep market-wide jumps\sep price jump\sep volatility jump\\[.25cm]
\end{keyword}

\end{abstract}

\end{frontmatter}
\thispagestyle{plain}

\section{Introduction}
\renewcommand*{\thefootnote}{\arabic{footnote}}
\setcounter{footnote}{0}
Understanding the relation between asset returns and volatility is among the most enduring and highly active research topics in finance. From an economic point of view, there seems to be a consensus that stock market returns and changes in volatility should be negatively correlated.\footnote{Some papers define the leverage effect as the correlation between returns and the {\em level} of volatility. \cite{d95} discusses the relation between the two definitions.}  
The linear, inverse return-volatility relationship is usually attributed to both changes in financial leverage and a time-varying risk premium; see \cite{b76},  \cite{fss87}, \cite{d95}, \cite{bw00} and \cite{blt06}. 
The financial leverage explanation motivates labeling the purely statistical relation between stock returns and volatility as the  ``leverage effect.'' 

Estimation of the leverage effect is challenging. \cite{afl13} document that the leverage effect fades out when using data sampled at increasing observation frequencies. In the framework of the Heston model, they show that discretization errors, volatility estimates and market microstructure noise bias the na\"{i}ve return-volatility correlation estimator towards zero. Recent research has tried hard to establish the leverage effect for intraday data. 

If the asset price and volatility processes have both Brownian and jump components, then the relation between returns and volatility splits into continuous and discontinuous parts. Continuous leverage refers to the relation between the Brownian components of the price and volatility processes. \cite{v12}, \cite{wm14}, \cite{aflwy16} and \cite{kx16} study measures of continuous leverage. These papers document a negative and usually time-varying continuous leverage effect. The discontinuous leverage effect (DLE) measures the relation between sizes of contemporaneous price and volatility jumps. \cite{br16} highlight the crucial importance of  both leverage components for asset pricing and risk management. Specifically, they show how a price-volatility cojump covariation affects return and variance risk premia. Their model estimates suggest that  discontinuous leverage explains about 25\% of the S\&P 500 return risk premium. This economically sizable proportion depends on the frequency of price-volatility cojumps, as well as the sign and magnitude of the covariation.

The existence of the DLE appears controversial, however. 
Several previous studies reached different conclusions regarding a DLE. 
\cite{jkm13} use truncated returns and increments of local spot volatility estimates to construct correlation statistics for one-minute S\&P 500 Exchange-Traded Funds (ETF) data from 2005 to 2011. These statistics indicate little evidence of a DLE. In contrast, \cite{br16} focus on a relatively small set of very large price jumps and a spot variance estimator based on infinitesimal cross-moments for high-frequency S\&P 500 futures from 1982 to 2009. Their parametric estimates suggest a strong DLE with correlations from -0.6 to -1. \cite{aflwy16} find that the DLE for five-second Dow Jones index data from 2003 to 2013 is usually different from zero. Their empirical analysis does not recover the sign and magnitude of the discontinuous leverage, however. Finally, \cite{tt11} use five-minute option-implied volatility index (VIX) data to evaluate volatility jumps in the S\&P 500 index from 2003 to 2008. The authors find that squared jumps in the S\&P 500 index are strongly positively correlated with jumps in the VIX. All these papers focus on stock market indexes, not individual stocks, and only use methods that are not robust to market microstructure noise.

Our paper makes both methodological and empirical contributions. We introduce novel methods to estimate and test the covariation of contemporaneous price and volatility jumps---denoted by \cite{aflwy16} as the DLE. A direct extension of our covariation estimator consistently estimates the corresponding correlation.  
\cite{aflwy16} derive a limit theorem for the DLE estimator that only applies to a setting without market microstructure noise. \cite{cop14} point out, however, that it is important to use noise-robust methods and thereby to avoid downsampling the data to lower observation frequencies. Downsampling may result in spurious jump detection and affect the accuracy of discontinuous leverage estimates. Using noise-robust estimators for jumps in log prices and volatility, we establish a stable central limit theorem under market microstructure noise for the DLE for finite activity price jumps or large jumps of an infinite activity jump component. We provide a consistent, asymptotic test for the presence of the DLE.

We estimate the covariation using only the physical measure, i.e., observed stock prices. DLE estimation requires three steps: price-jump localization, price-jump estimation and estimation of changes in the spot volatility process at price-jump times. Under noise, none of the three steps is standard. We use spectral methods in all three steps. \cite{reiss} introduces spectral estimation of the quadratic variation from noisy observations. \cite{bhmr14} and \cite{ab15} establish the asymptotic efficiency of spectral estimators of the integrated volatility matrix in the multivariate case with noisy and non-synchronous observations. \cite{bhmr17} propose a related spot volatility estimator. Although spectral and the popular pre-average estimators have some similarities, they belong to different classes of estimators; see Remark 1 of \cite{bw16}. Our theoretical contribution is to provide methods to detect and estimate price jumps and to combine the three steps to infer the DLE. 

To detect price-jump times, we refine the adaptive thresholding approach of \cite{bw15}. We construct an argmax-estimator, such as is often used in change-point analysis. This refinement of the jump localization is motivated by the fact that estimation of price jumps becomes more difficult in cases where the jump times are not precisely determined. See \cite{v14} for a related problem. To estimate the price-jump size at a detected jump time, we first review the pre-average method of \cite{lm12} that extends the \cite{lm8} approach to a model with market microstructure noise. While \cite{lm12} mainly focus on a global test for jumps, we focus on local jump estimates. We generalize their stable central limit theorem from a jump diffusion to more general semimartingale models. Estimating the entire quadratic variation with jumps or testing for jumps over a whole day are related yet different problems. \cite{jpv10} and \cite{Koike2017} have developed rate-optimal consistent pre-average estimators for the quadratic variation and \cite{bw15} provide spectral estimators for this purpose. While these methods do not recover individual price jumps, the \cite{lm12} method utilizes natural \emph{local} average statistics to address inference on price jumps under noise. The pre-average method attains the optimal rate of convergence for local price-jump estimation. 
As one ingredient of the price-jump localization, we exploit the simple structure of these pre-average statistics that permits an asymptotic theory based on Gaussian approximations. Using spectral local statistics for price-jump estimation, we derive a superior estimator with a smaller variance than the pre-average estimator. The asymptotic variance of the spectral estimator attains the asymptotic lower bound. Thus, we provide the first feasible, asymptotically efficient estimator of price jumps from noisy observations.
To estimate changes in the spot volatility at a price-jump time, we employ the jump-robust techniques of Bibinger and Winkelmann (2016). Finally, we plug the price-jump and volatility-jump estimates into the DLE statistic of \cite{aflwy16}. 

Our methods provide new empirical evidence about the DLE for 320 individual stocks, which were actively traded at the NASDAQ stock exchange from 2010 to 2015. We find no prevalent evidence of an unconditional DLE in individual stock data, but we identify two forces that prevent significant unconditional discontinuous leverage estimates: First, while downward price jumps are usually negatively correlated with contemporaneous volatility jumps, upward price jumps are positively correlated with contemporaneous volatility jumps. Second, market jumps, i.e.\ price jumps that coincide with jumps of a market portfolio, display a strong DLE. In contrast, idiosyncratic price jumps, which occur without a contemporaneous jump of the market portfolio, are associated with a much weaker DLE. We establish an economically and statistically significant relation between cojumps in stock prices and volatility by conditioning on the sign of price jumps and whether those jumps are systematic or idiosyncratic. Apart from the sign of the DLE, we investigate its magnitude in cross-sectional regressions. We show that firms' debt-to-equity ratios do not explain much of the cross-section of DLE estimates or the correlations of price-volatility cojumps. In contrast, the volatility levels of individual firms strongly explain the magnitude of DLE estimates.

Our failure to find an unconditionally negative DLE is consistent with the asset pricing models of \cite{pv12} in which specific events trigger jumps. That is, the continuous leverage effect and the DLE are fundamentally different in that model. Their learning model implies that changes in monetary or government policy trigger market-wide price and volatility cojumps, where the uncertainty about the impact of a new policy regime on the profitability of private firms always raises volatility, regardless of the effect on prices. News that causes asset prices to jump up while causing volatility to jump down is incompatible with their model. Our results are also consistent with \cite{p17}, who studies systematic and nonsystematic risk factors in S\&P500 high-frequency firm data. That is, we confirm that the DLE appears predominantly for systematic risk, while being weaker and more often nonsignificant for idiosyncratic risk. 
 
The rest of the paper is organized as follows. Section 2 introduces the model and assumptions. Section 3 presents the price-jump estimators, spot volatility estimation and the DLE estimator. We compare the spectral approach for price jumps with the \cite{lm12} pre-average estimator. Section 4 provides Monte Carlo evidence and Section 5 the empirical findings. Section 6 concludes. The Appendix contains the proofs.

\section{Statistical model and assumptions}\label{sec2}
We work with a very general class of continuous-time processes, namely It\^o semimartingales. Its implicit no-arbitrage properties make it the most popular model for log-price processes in financial econometrics. The model is formulated for a log price, $X_t$, and its volatility, $\sigma_t$, over a fixed time period $t\in[0,1]$, on some filtered probability space $(\Omega,\mathcal{F},(\mathcal{F}_t),\mathbb{P})$:   
\begin{eqnarray}X_t&=&X_0+\int_0^t b_s\,ds+\int_0^t\sigma_s\,dW_s + \int_0^t \int_{\mathbb{R}}\delta(s,z)\1_{\{|\delta(s,z)|\leq 1 \}}(\mu-\nu)(ds,dz) \nonumber \\ && +\int_0^t \int_{\mathbb{R}}\delta(s,z)\1_{\{|\delta(s,z)|> 1 \}} \mu(ds,dz)\,,\label{sm} \end{eqnarray}
with a standard Brownian motion $(W_s)$, the jump size function $\delta$, defined on  $\Omega\times \mathbb{R}_+\times \mathbb{R}$, and the Poisson random measure $\mu$, which is compensated by $\nu(ds,dz)=\lambda(dz)\otimes ds$ with a $\sigma$-finite measure $\lambda$. We write $\Delta X_t=X_t-X_{t-}$ with $X_{t-}=\lim_{s< t,s\rightarrow t} X_s$ for the process of jumps in $(X_t)$ and $\Delta\sigma_t^2=\sigma^2_t-\sigma^2_{t-}$ for jumps of the squared volatility. Our notation follows that of \cite{jp12}. We impose mild regularity assumptions on the characteristics of $X_t$.
\begin{assump}\label{sigma}\onehalfspacing
The drift $(b_t)_{t\ge 0}$  is a locally bounded process. The volatility never vanishes, $\inf_{t\in[0,1]}\sigma_t>0$ almost surely. 
For all $0\leq t+s\leq1$, $t\ge 0$, some constants $C_n,K_n>0$, some $\alpha>1/2$ and for a sequence of stopping times $T_n$, increasing to $\infty$, we have that 
\begin{align}
\label{JM2}\Big|\E\big[\sigma_{(t+s)\wedge T_n}-\sigma_{t\wedge T_n}\,|\mathcal{F}_t\big]\Big| & \le C_n\,s^{\alpha}\,,\\
\label{JM}\E\Big[\sup_{\mathfrak{t}\in[0,s]}|\sigma_{(\mathfrak{t}+t)\wedge T_n}-\sigma_{t \wedge T_n}|^2\Big] & \le K_n\,s\,.\end{align}
\end{assump}
Assumption \ref{sigma} requires some smoothness of the volatility process. It does not exclude volatility jumps, only fixed times of discontinuity are excluded. We impose the following regularity condition on the jumps.
\begin{assump}\label{jumps}\onehalfspacing
Assume for the predictable function $\delta$ in \eqref{sm} that \(\sup_{\omega,x}|\delta(t,x)|/\gamma(x)\) is locally bounded with a non-negative, deterministic function $\gamma$ that satisfies
\begin{align}\label{BG}\int_{\mathbb{R}}(\gamma^r(x)\wedge 1)\lambda(dx)<\infty\,.\end{align}
\end{assump}
The index $r,\, 0\le r\le 2$, in \eqref{BG} measures the jump activity. Smaller values of $r$ make Assumption \ref{jumps} more restrictive. In particular, $r=0$ results in finite-activity jumps and $r=1$ implies that jumps are summable.
\begin{remark}\onehalfspacing
Assumption \ref{sigma} is satisfied in a very general model, where the volatility process $\sigma_t$ is an It\^{o} semimartingale 
\begin{eqnarray}
\sigma_t&=&\sigma_0+\int_0^t \tilde b_s\,ds+\int_0^t\tilde \sigma_s\,d\tilde{W}_s+ \int_0^t \int_{\mathbb{R}}\tilde \delta(s,z)\1_{\{|\tilde \delta(s,z)|\leq 1 \}}(\mu-\nu)(ds,dz) \nonumber\\ &&+\int_0^t \int_{\mathbb{R}}\tilde \delta(s,z)\1_{\{|\tilde \delta(s,z)|> 1 \}}\mu(ds,dz)\,\label{smvola},
\end{eqnarray}
with a standard Brownian motion $(\tilde W_s)$, when the characteristics in \eqref{smvola} are locally bounded and when an analogous condition as \eqref{BG} holds for $\tilde\delta $ in \eqref{smvola} with $r=2$. We may use the same $\mu$ in \eqref{sm} and \eqref{smvola}, such that $\mu$ is a jump measure on $\mathbb{R_+}\times \mathbb{R}$ governing the jumps in the log price and its volatility. The predictable functions, $\delta$ and $\tilde \delta$, defined on $\Omega\times \mathbb{R}_+\times \mathbb{R}$, then determine common jumps of $\sigma_t$ and $X_t$. Whenever $\delta\tilde\delta\equiv 0$, there is no price-volatility cojump. Our asymptotic theory and Assumption \ref{sigma} allow for generalizations of \eqref{smvola}. For instance, long-memory fractional volatility components can be included. Thus, our theoretical setup includes almost any popular stochastic volatility model that allows for both continuous and discontinuous leverage effects. 
\end{remark}
In practice, one cannot observe the efficient price \eqref{sm} directly and one must account for market microstructure noise in analyzing price and volatility jumps. To efficiently exploit available high-frequency prices, we posit a latent discrete observation model with noise:
\begin{align} \label{obs}\text{Observe} ~~~{Y_{t_i^n}},i=0,\ldots,n,~~\text{with}~~Y_t=X_t+\epsilon_t\,,~~~~~~~~~~\end{align}
where $\epsilon_t$ captures the market microstructure noise. 
We use the typical notation, $\Delta_i^n Y={Y_{t_i^n}-Y_{t_{i-1}^n}},i=1,\ldots,n$, for noisy returns and analogous notation for the processes $(X_t)$ and $(\epsilon_t)$. In our baseline setup, market microstructure noise is a white noise process $(\epsilon_t)_{t\ge 0}$, independent of $X_t$, with $\E[\epsilon_t]=0$ and $\E[\epsilon_t^2]=\eta^2$, as well as $\E[\epsilon_t^{4+\delta}]<\infty$ for some $\delta>0$, for all $t\in[0,1]$. The process $Y_t$ is accommodated on the product space $(\bar\Omega,\mathcal{G},(\mathcal{G}_t),\bar\P)$, where $\mathcal{G}_t=\mathcal{F}_t\otimes\sigma(\epsilon_s,s\le t)$ contains information about the signal and noise. Below we extend the model to more general setups with serially correlated, heteroscedastic noise. Because we apply our methods to locally infer price and volatility jumps of individual stock prices, non-synchronicity of the multivariate data is of less importance here.

 \section{Inference on the discontinuous leverage effect\label{sec:3}}
The DLE is defined as the covariation of contemporaneous price and volatility jumps. We estimate it in three steps. We first address noise-robust estimation of price jumps in Section \ref{sec:3.1}, then  we turn to noise-robust estimation of spot-volatility changes in Section \ref{sec:3.2}. In Section \ref{sec:3.3}, we show how to detect a priori unknown jump times in noisy data and how to refine price-jump estimation for DLE estimation in this case. The covariation of the price-jump and spot-volatility estimates at detected jump times gives the estimated DLE.

 \subsection{Price-jump estimation\label{sec:3.1}}
 \subsubsection{Local jump estimator and test using pre-averaged log prices}
 Consider  the statistic
\begin{align}\label{lm}T^{LM}(\tau;\Delta_1^n Y,\ldots,\Delta_n^n Y)=\hat P(t_l^n)-\hat P(t_{l-M_n}^n)\,,\,l=\lfloor \tau n\rfloor+1\,,\end{align}
at a (stopping) time $\tau\in(0,1)$ and with pre-processed price estimates 
\begin{align}\label{locest}\hat P(t_j^n)=M_n^{-1}\sum_{i=j}^{(j+M_n-1)\vee n}Y_{t_i^n}\,.\end{align}
\cite{lm12} propose a test for price jumps at time $\tau$ based on \eqref{lm}. The window length for the pre-averaging is $M_n=c\sqrt{n}$ with a proportionality constant $c$. The following proposition generalizes Lemma 1 in \cite{lm12}, where the authors assume that they observe discrete, noisy observations from a jump-diffusion model.
\begin{prop}\label{proplm}\onehalfspacing Under Assumption \ref{sigma} and Assumption \ref{jumps} with $r<4/3$ for equidistant observations, $t_i^n=i/n$, the Lee-Mykland statistic \eqref{lm} obeys the stable\footnote{Stable means stable convergence in law with respect to $\mathcal{F}$.} central limit theorem,
\begin{align}\label{stablecltlm}\sqrt{M_n}\,\big(T^{LM}(\tau;\Delta_1^n Y,\ldots,\Delta_n^n Y)-\Delta X_{\tau}\big)\stackrel{(st)}{\longrightarrow} MN\Big(0,{\frac13 (\sigma_{\tau}^2+\sigma_{\tau-}^2)}\,c^2+2\eta^2\Big)\,,\end{align}
as $n\rightarrow\infty$, where $MN$ stands for mixed normal. 
\end{prop}
Thus, in case of a price jump at $\tau$, \eqref{lm} consistently estimates the price-jump size. The central limit theorem accounts for a contemporaneous volatility jump. If there is no volatility jump, then $\sigma_{\tau}^2=\sigma_{\tau-}^2$ in \eqref{stablecltlm}. With the null hypothesis, $ \Delta X_{\tau}=0$, and alternative, $|\Delta X_{\tau}|>0$, Proposition \ref{proplm} facilitates a consistent test for a jump in the stock price at time point $\tau\in (0,1)$. 

\subsubsection{Local jump estimator and test using spectral statistics} 
To estimate price jumps using spectral statistics, we consider an orthogonal system of sine functions that are localized on a window around $\tau$:
\begin{align}\label{spectral}\Phi_{j,\tau}(t)=\sqrt{\frac{2}{h_n}}\sin{\big(j\pi h_n^{-1}(t-(\tau-h_n/2))\big)}\mathbbm{1}_{[\tau-h_n/2,\tau+h_n/2]}(t)\,,\;j\ge 1\,.\end{align}
Asymptotically efficient volatility estimation from noisy observations \eqref{obs} motivates consideration of local averages of noisy log prices in the frequency domain; see \cite{reiss} and \cite{bhmr14}. Intuitively, \emph{spectral statistics},
\begin{align}\label{spec}S_j(\tau)=\sum_{i=1}^n\Delta_i^n Y\Phi_{j,\tau}((t_{i-1}^n+t_i^n)/2)\,,\;j\ge 1\,,\end{align}
maximize the local information load about the signal process and thereby allow for local estimates of the efficient prices: $X_{\tau}$ and $X_{\tau-}$. The scaling factor in front of the sine in \eqref{spectral} ensures that $\int_{\tau-h_n/2}^{\tau+h_n/2}\Phi_{j,\tau}^2(t)\,dt=1$. 
We propose the following statistic: 
\begin{align}\label{bw}\mathcal{T}(\tau;\Delta_1^n Y,\ldots,\Delta_n^n Y)=\sum_{j=1}^{J_n}(-1)^{j+1}a_{2j-1}S_{2j-1}(\tau)\sqrt{h_n/2}\,,\end{align}
with weights $(a_{2j-1})_{j\ge 1}$, to infer price jumps. \eqref{bw} is a rescaled weighted sum of spectral statistics over odd spectral frequencies up to some spectral cut-off frequency $2J_n-1$. Excluding even frequencies and alternating the signs of addends facilitate a consistent estimation of price jumps $\Delta X_{\tau}$, as in \eqref{stablecltlm} above. 

The window length is set to be $h_n=\kappa \log{(n)}/\sqrt{n}$ for some constant $\kappa$. Despite the logarithmic factor, the window length resembles the one in \eqref{locest}. We derive optimal oracle weights by minimizing the variance, which depend on time through the volatility $\sigma_t$. Yet, under Assumption \ref{sigma}, the error of approximating $(\sigma_t^2)$ constant on $[\tau-h_n/2,\tau)$ and $[\tau, \tau+h_n/2]$ is asymptotically negligible. Then, as in the weighted least squares approach, this leads to optimal weights
\[a_j\propto 1/\var(S_{j}(\tau))\,.\]
In order to consistently estimate the jump $(X_{\tau}-X_{\tau-})$, we set $\sum_{j=1}^{J_n}a_{2j-1}=1$ such that
\begin{align}\label{weights}a_{2j-1}&=\frac{\big(\var\big(S_{2j-1}(\tau)\big)^{-1}}{\Big(\sum_{u=1}^{J_n}\big(\var\big(S_{2u-1}(\tau)\big)\big)^{-1}\Big)}\\
&\notag =\frac{(\frac12(\sigma_{\tau}^2+\sigma_{\tau-}^2)+\pi^2(2j-1)^2h_n^{-2}n^{-1}\eta^2)^{-1}}{\Big(\sum_{u=1}^{J_n}(\frac12(\sigma_{\tau}^2+\sigma_{\tau-}^2)+\pi^2(2u-1)^2h_n^{-2}n^{-1}\eta^2)^{-1}\Big)}\,.\end{align}
For an adaptive method, we estimate these oracle optimal weights by plugging in the estimated noise variance, 
\begin{subequations}
\begin{align}\label{pre1}\hat\eta^2=-n^{-1}\sum_{i=1}^{n-1}\Delta_i^n Y\Delta_{i-1}^n Y=\eta^2+\mathcal{O}_{\bar\P}\big(n^{-1/2}\big)\,,\end{align}
and the pre-estimated spot squared volatility,
\begin{align}\hat\sigma_{\tau-,pil}^2&=\frac{r_n}{J_p}\sum_{k=1}^{r_n^{-1}}\sum_{j=1}^{J_p}\big(S_j^2(\tau-kh_n)-\pi^2j^2h_n^{-2}n^{-1}\hat\eta^2\big) \nonumber\\
&\hspace*{1cm}\notag\times \1\Big(\Big|(J_p)^{-1}\sum_{j=1}^{J_p}\big(S_j^2(\tau-kh_n)-\pi^2j^2h_n^{-2}n^{-1}\hat\eta^2\big)\Big|\le u_n\Big) \\ &=\sigma_{\tau-}^2+\mathcal{O}_{\bar\P}\big(n^{-1/8}\big)\,,\label{pre2}\end{align} 
\end{subequations} 
with $r_n^{-1}=\mathcal{R}n^{1/4}$ for a constant $\mathcal{R}$, a threshold sequence $u_n=h_n^{\varpi},\, 0<\varpi<1$, and maximal spectral frequency, $J_p$, leading to the above rate-optimal estimators under Assumptions \ref{sigma} and \ref{jumps} with $r<3/2$. The notation $S_j^2(\tau-kh_n)$ refers to squared spectral statistics computed from $r_n^{-1}$ bins with sine functions centered around times, $\tau-kh_n$, before $\tau$. $\hat\sigma_{\tau,pil}^2$ is the analog of \eqref{pre2}, replacing $\tau-kh_n$ by $\tau+kh_n$. \cite{bw16} detail the construction and prove the asymptotic properties of pre-estimators \eqref{pre1} and \eqref{pre2} and also suggest how to choose $\mathcal{R}$ and $J_p$.

\begin{figure}[t]
\begin{center}
\includegraphics[width=9cm]{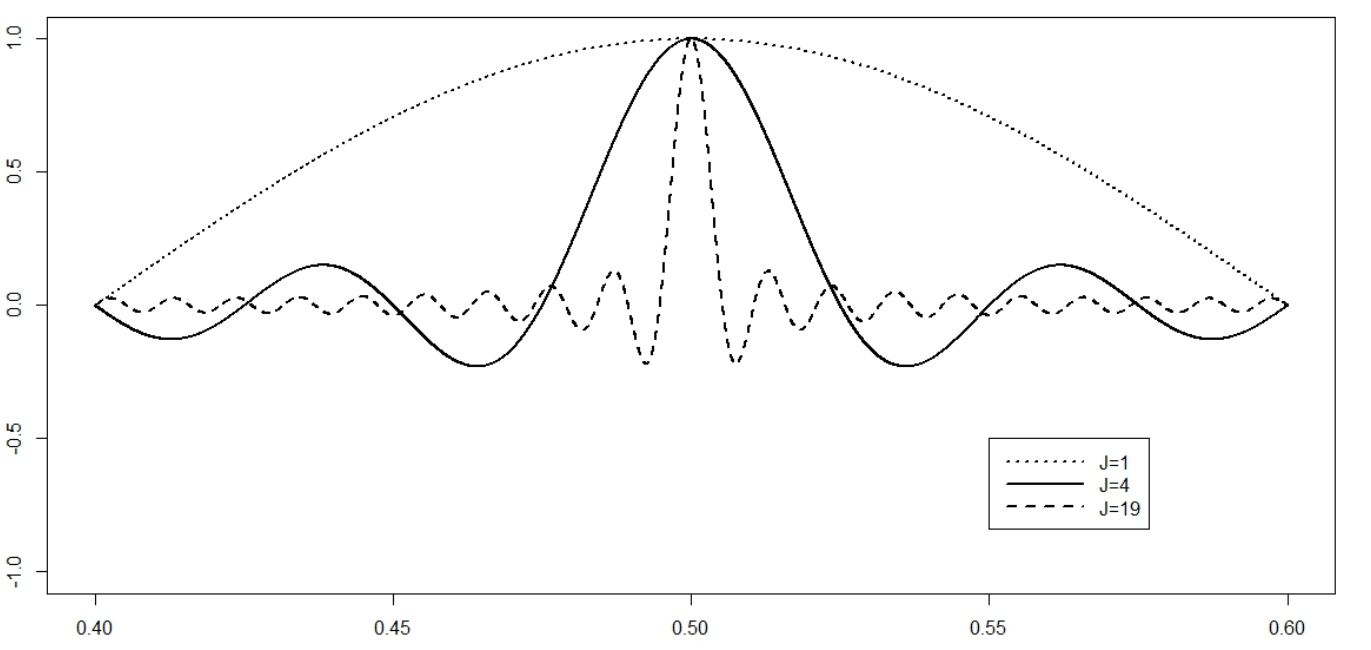}\\
\includegraphics[width=9cm]{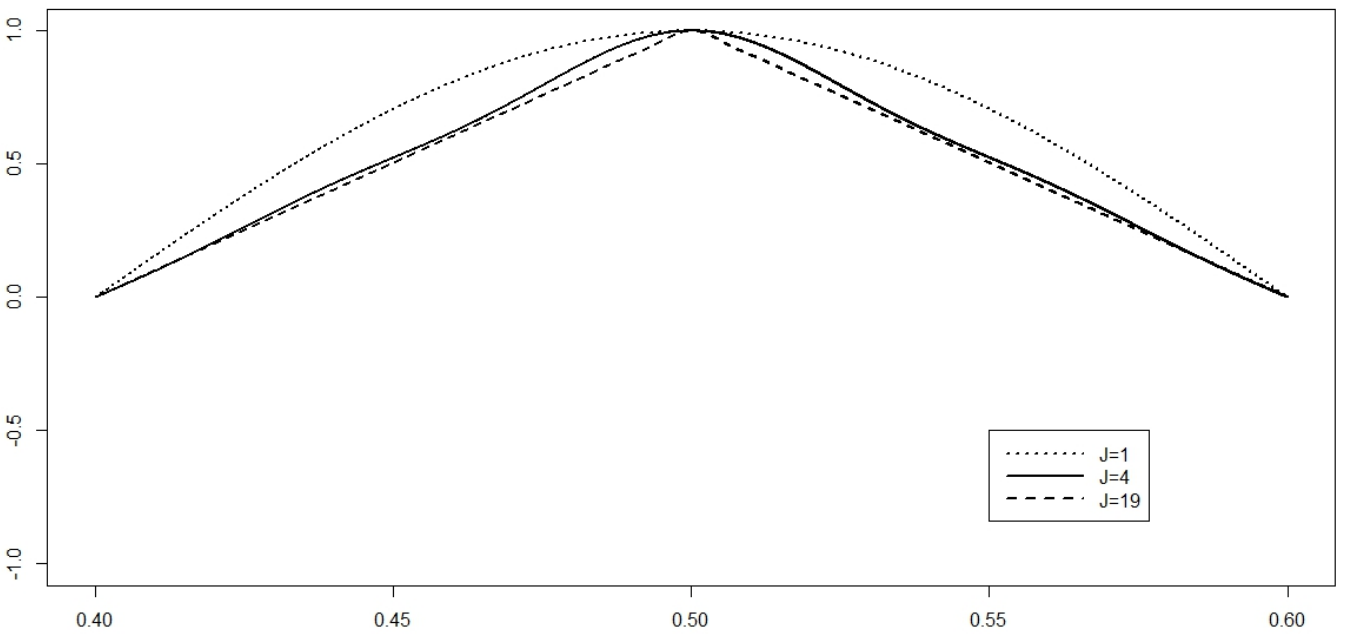}
\end{center}
\caption{\label{Fig1} Example for $\sum_{j=1}^JJ^{-1}(-1)^{j+1}\Phi_{2j-1,\tau}(t)\sqrt{h_n/2}$ (top) and, with oracle optimal weights $a_{2j-1}$ from the simulation setup in Section \ref{sec:sim}, $\sum_{j=1}^Ja_{2j-1}(-1)^{j+1}\Phi_{2j-1,\tau}(t)\sqrt{h_n/2}$ (bottom) as functions of time $t$ for three values of $J$ on window $[\tau-h_n/2,\tau+h_n/2], \tau=0.5,h_n=0.2$. The plots demonstrate that if a jump occurs on the interval, and not sufficiently close to $\tau$, the estimation of $\Delta X_{\tau}$ becomes biased for the actual jump.}
\end{figure}

Next, we state asymptotic results for $\mathcal{T}(\tau;\Delta_1^n Y,\ldots,\Delta_n^n Y)$, which refers to statistic \eqref{bw} with estimated optimal weights.
\begin{prop}\label{propbw}\onehalfspacing
Under Assumption \ref{sigma} and Assumption \ref{jumps} with $r<4/3$ for equidistant observations, $t_i^n=i/n$, our statistic \eqref{bw} obeys the stable central limit theorem as $n\rightarrow\infty$ and $J_n\rightarrow\infty$: 
\begin{align}\label{stablecltbw}n^{1/4}\,\big(\mathcal{T}(\tau;\Delta_1^n Y,\ldots,\Delta_n^n Y)-\Delta X_{\tau}\big)\stackrel{(st)}{\longrightarrow} MN\Big(0,2\Big(\frac{\sigma_{\tau}^2+\sigma_{\tau-}^2}{2}\Big)^{1/2}\eta\Big)\,.\end{align}
\end{prop}
In the case of no volatility jump at $\tau$, $\sigma_{\tau}=\sigma_{\tau-}$ and the asymptotic variance is $2\sigma_{\tau}\eta$.
Finally, we extend Proposition \ref{propbw} to a more realistic model that incorporates serially correlated, heteroscedastic noise and non-regular sampling.
\begin{assump}\label{eta}\onehalfspacing
Assume the existence of a differentiable, cumulative distribution function $F$ that determines the observation times via a quantile transformation, $t_i^n=F^{-1}(i/n),i=0,\ldots,n$. Assume $(F^{-1})'$ is $\alpha$-H\"older continuous for some $\alpha>1/2$, i.e., $|(F^{-1})'(t)-(F^{-1})'(s)|\le |t-s|^{\alpha}$ for all $s,t$.\\
The noise process $(\epsilon_{t})$ is independent of $X$. For all $t$, we have $\E[\epsilon_{t}]=0$ and $\E\big[\epsilon_t^{4+\delta}\big]<\infty$, for some $\delta>0$. Further, assume $(\epsilon_{t_i^n})$ is an $R$-dependent process, such that $\cov(\epsilon_{t_i^n},\epsilon_{t_{i+u}^n})=0$ for $u>R$ and some $R<\infty$, then 
the long-run variance process converges as follows:
\begin{align}\label{lrvar}\sum_{l=-\lfloor t n\rfloor }^{n-\lfloor t n\rfloor}\cov\big(\epsilon_{\lfloor t n\rfloor },\epsilon_{\lfloor t n\rfloor +l}\big)\rightarrow \eta^2_t~,\end{align}
for $t\in[0,1]$, uniformly in probability. The process $(\eta_t^2)_{t\in[0,1]}$ is locally bounded and satisfies, for all $t,(t+s)\in[0,1]$, the mild smoothness condition:
\begin{align} \label{etasmooth}|\eta_{t+s}^2-\eta^2_t|\le K s^{\alpha}\,.\end{align}
The noise does not vanish: $\eta_t^2>0$ for all $t\in[0,1]$. 
\end{assump}
\begin{prop}\label{propbwgen}\onehalfspacing
Under Assumptions \ref{sigma}, \ref{jumps} with $r<4/3$ and \ref{eta}, the statistic \eqref{bw} obeys the stable central limit theorem as $n\rightarrow\infty$ and $J_n\rightarrow\infty$ 
\begin{align}\label{stablecltbwgen}n^{1/4}\,\big(\mathcal{T}(\tau;\Delta_1^n Y,\ldots,\Delta_n^n Y)-\Delta X_{\tau}\big)\stackrel{(st)}{\longrightarrow} MN\Big(0,2\Big(\frac{\sigma_{\tau}^2+\sigma_{\tau-}^2}{2}\Big)^{1/2}\hspace*{-.025cm}\eta_{\tau} \big((F^{-1})'(\tau)\big)^{1/2}\Big).\end{align}
\end{prop}
In analogy to Proposition \ref{proplm}, Propositions \ref{propbw} and \ref{propbwgen} show the consistency of the spectral jump estimator and give a consistent test for a price jump at time $\tau$. One can construct standardized, feasible versions of \eqref{stablecltbw} and \eqref{stablecltbwgen} by inserting spot squared volatility and long-run noise variance estimators. See \eqref{pre1} and \cite{bhmr17} for such estimators. In fact, the pre-estimation of optimal weights also provides estimates of the variances of \eqref{bw}.
The asymptotic variance of the Lee-Mykland statistic in \eqref{stablecltlm} generalizes  to $(1/3)(\sigma_{\tau}^2+\sigma_{\tau-}^2)(F^{-1})'(\tau) c^2+2\eta_{\tau}^2$ under the conditions from Proposition \ref{propbwgen}. \cite{lm12} provide a generalization to $R$-dependent noise using sub-sampling, and this directly applies to our general setup with Assumption~\ref{eta}. The spectral price-jump estimator \eqref{bw} and the pre-average jump estimator \eqref{lm} have the same optimal convergence rate and similar asymptotic properties.
\begin{remark}\onehalfspacing Writing \eqref{stablecltlm} with rate $n^{1/4}$ instead of $M_n^{1/2}$, the variance of the Lee-Mykland estimator \eqref{lm} with $\sigma_{\tau}=\sigma_{\tau-}$ becomes $\frac23\sigma^2_{\tau}c+2\eta^2c^{-1}$. The variance is minimized by the constant $c=\sqrt{3}\eta\sigma_\tau^{-1}$, which yields $4\sigma_\tau\eta/\sqrt{3}$ in \eqref{stablecltlm}. Since $4/\sqrt{3}\approx 2.31$, this optimized variance of (an infeasible) Lee-Mykland estimator is about 16\% larger than the variance $2\sigma_{\tau}\eta$ of the spectral estimator in \eqref{stablecltbw}. Moreover, according to the LAN result of \cite{Koike2017}, the latter is optimal. That is, the variance of the spectral estimator coincides with a lower bound for the asymptotic variance, which is given by the inverse of the Fisher information from Proposition 5.2 of \cite{Koike2017}. Our estimator is hence the first feasible, asymptotically efficient estimator for price jumps in the semimartingale model with market microstructure noise.
\end{remark} 
We caution, however, that estimates via spectral statistics \eqref{bw} and pre-averages \eqref{lm} are biased when a jump is not located close to time $\tau$ but instead close to the edges of the local window. Figure \ref{Fig1} illustrates this. The bias for the Lee-Mykland estimator is linear. This effect directly relates to the so-called ``pulverisation'' of jumps by pre-averages described in \cite{mykzhang16}. For our statistic, the bias hinges on the weights and the spectral cut-off.  The lower panel of Figure \ref{Fig1} reveals that the bias is similar for both methods. The bias becomes important when studying price jumps at a priori unknown times, such as when one is estimating the DLE.  Section \ref{sec:3.3} discusses our solution. Related problems by not knowing the exact timing of jumps arise and have been addressed in different ways in \cite{v14} and \cite{bw15}.
\subsection{Spot volatility estimation\label{sec:3.2}} \label{pjvjestimats}
We estimate the contemporaneous volatility adjustment to a price jump at time $\tau\in(0,1)$.
We employ the spectral spot squared volatility estimators of \cite{bw16}, smoothed over local windows before $\tau$  and after $\tau$, to consistently estimate the volatilities $\sigma^2_{\tau}$ and $ \sigma^2_{\tau-}$. Based on estimates of the oracle optimal weights
\begin{align}\label{orweights}w_{jk}=I_k^{-1}I_{jk}=\frac{\Big(\sigma^2_{(k-1)h_n}+\pi^2j^2h_n^{-2}\tfrac{\eta_{(k-1)h_n}^2}{n}\Big)^{-2}}{\sum_{m=1}^{J_n}\Big(\sigma^2_{(k-1)h_n}+\pi^2m^2h_n^{-2}\tfrac{\eta_{(k-1)h_n}^2}{n}\Big)^{-2}}\,,\end{align}
inserting spot squared volatility and noise variance estimators, with
\begin{align}\label{zetaad}\zeta_k^{ad}(Y)=\sum_{j=1}^{J_n} \hat w_{jk}\Big(S_{jk}^2-\pi^2j^2h_n^{-2}\tfrac{\hat\eta_{(k-1)h_n}^2}{n}\Big)\,,\end{align}
the spectral estimator of the spot squared volatility at time $\tau-$ is
\begin{align}\label{lspot}
\hat \sigma^2_{\tau-}=r_n\sum_{k=\lfloor sh_n^{-1}\rfloor -r_n^{-1}}^{\lfloor sh_n^{-1}\rfloor-1}\zeta_k^{ad}(Y)\1_{\{h_n|\zeta_k^{ad}(Y)|\le u_n\}}\,.\end{align}
To estimate the noise variance $\eta^2$ and pre-estimate the spot squared volatility in \eqref{orweights}, we use  \eqref{pre1} and \eqref{pre2}, respectively.  To obtain \eqref{zetaad}, we adapt $S_{jk}=S_j((k-1/2)h_n)$ from \eqref{spec}. Analogously to $\hat\sigma_{\tau-}^2$, $\hat\sigma_{\tau}^2$ is defined by summing over $k\in\{\lfloor sh_n^{-1}\rfloor+1,\ldots,\lfloor sh_n^{-1}\rfloor +r_n^{-1}\}$. The theory by \cite{bw16} renders the following result:
\begin{cor}\label{testvolajump}\onehalfspacing Under Assumptions \ref{sigma} and \ref{jumps} with $r<3/2$ for equidistant observations, $t_i^n=i/n$, and under Assumption \ref{eta} for the noise, the statistics \eqref{lspot} with $r_n\propto n^{-\beta}\log(n)$ satisfy
\begin{align}\label{testlt}n^{\beta/2}\,\big(\big(\hat \sigma^2_{\tau}-\hat \sigma^2_{\tau-}\big)-\Delta\sigma^2_{\tau}\big)\stackrel{(st)}{\longrightarrow}MN\big(0,8(\sigma_s^{3}+\sigma_{s-}^{3})\eta_s\big)\,\end{align}
for all
\begin{align}\label{beta}0<\beta<\left(1/4\,\wedge\,\varpi\Big(1-\frac{r}{2}\Big)\right)\,, \end{align}
with $\varpi$ from the truncation sequence $u_n$, such that we come arbitrarily close to the optimal rate $n^{1/8}$ in \eqref{testlt}.  
\end{cor}
Theorem 10.30 of \cite{sahaliajacod} provides a related result for volatility jump estimation without microstructure noise, where the first line in their equation (10.81) for one fixed point in time corresponds to our result under condition \eqref{beta}. The corollary provides an asymptotic test of the hypothesis of no volatility jump, $\Delta\sigma_{\tau}^2=0$, against the alternative that $\Delta\sigma_{\tau}^2\ne 0$. The statistic (27) in \cite{bw16} gives an efficient test. For non-equidistant observations, the noise level in \eqref{testlt} includes $((F^{-1})'(\tau))^{1/2}$, analogous to \eqref{stablecltbwgen} for price jumps.

 \subsection{Discontinuous leverage effect\label{sec:3.3}}
This section introduces a covariation measure for contemporaneous price and volatility jumps that combines the above spectral jump and volatility estimators. Our covariation measure is related to that of \cite{aflwy16} who introduce the following as the \emph{tail discontinuous leverage effect} in their equation (2.7):
\begin{align}\label{dle}[X,\sigma^2]^d_T(a)=\sum_{s\le T}\Delta X_{s}\big(\sigma^2_{s}-\sigma^2_{s-}\big)\,\1_{\{|\Delta X_s|>a\}}\,.\end{align}
Based on our local methods and setting $T=1$, we consider the DLE estimator 
\begin{align}
\widehat{[X,\sigma^2]}^d_1(a)=\sum_{k=2}^{h_n^{-1}-1}\widehat{\Delta X}_{\hat\tau_k}\big(\hat\sigma^2_{\hat\tau_k}-\hat\sigma^2_{\hat\tau_k-}\big)\,\1_{\{\Delta_k\widehat{[X,X]}>a^2\vee u_n\}}\,,\label{lev}\end{align}
where $\widehat{\Delta X}_{\hat\tau_k}$ is the estimated log-price jump \eqref{bw} and $\hat\sigma^2_{\hat\tau_k-}$ and $\hat\sigma^2_{\hat\tau_k}$ are the spot volatility estimates \eqref{lspot}. Only finitely many addends with (large) price jumps in \eqref{lev} are non-zero. \eqref{lev} makes it apparent that we need to detect these unknown price-jump times to estimate the DLE. With a fixed $a> 0$, or with $a=0$ in case of finite activity jumps, we first use a thresholding procedure to locate bins $((k-1)h_n,kh_n)$ that contain a (large) price jump. We apply a bin-wise threshold, $u_n(kh_n)=2\log (h_n^{-1})h_n\hat\sigma_{(k-1)h_n,pil}^2$, with pre-estimated squared volatility, as defined in \eqref{pre2}. The moving threshold accounts for intraday volatility patterns. To estimate changes in the quadratic variation, $ \Delta_k\widehat{[X,X]}$, on bins with a price jump, we adapt the statistics from Section 3.1.3 of \cite{bw15} and define
\begin{align}\label{zetatilde}\widetilde\zeta^{ad}_{k,l}&=\sum_{j\in\mathcal{J}_n} \hat w_{jk}\Big(\tfrac12 S_{jk}^2+\tfrac12\tilde S_{jl}^2-\pi^2j^2h_n^{-2}\tfrac{\hat\eta_{(k-1)h_n}^2}{n}\Big)\,, \nonumber \\ \widetilde\zeta^{ad}_{k}&=\max\big(\widetilde\zeta^{ad}_{k,k},\widetilde\zeta^{ad}_{k,k+1}\big)\,,\end{align} 
by summing over the set $\mathcal{J}_n$ of odd numbers up to the cut-off $J_n$ and with spectral statistics $\tilde S_{jk}=S_j((k-1)h_n)$ shifted by $h_n/2$ in comparison to $S_{jk}=S_j((k-1/2)h_n)$. This adjustment of \eqref{zetaad} allows for unbiased estimation of the increase in the jump variation on bins with jumps. Due to the overlapping nature of shifted bins and the maximum operator in \eqref{zetatilde}, a jump on a bin also affects a neighboring bin. The weighting of a jump on a neighboring bin is always smaller, however, than the weighting on the bin containing the jump. Thus, the increment in jump variation on a bin containing a jump is estimated by
\begin{align}\label{quadad}\Delta_k\widehat{[X,X]}=h_n\,\widetilde\zeta^{ad}_{k}\,\1_{\{\widetilde\zeta^{ad}_{k}>\max (\widetilde\zeta^{ad}_{k-1},\widetilde\zeta^{ad}_{k+1})\}}\,.\end{align}
The thresholding procedure detects asymptotically small bins with jumps. However, following \cite{bw15}, these bin-widths decay with order $n^{-1/2}$ and are of the same sizes as the bins in which the price-jump estimation is conducted. To solve the bias problem in price-jump estimation at unknown times, we determine a price-jump time $\tau$ more precisely and discuss how to apply our price-jump statistic, \eqref{bw}. While we can directly estimate volatility jumps from \eqref{lspot}, we need to adjust the price-jump estimation to obtain an overall consistent estimator of the DLE.
 

To determine the jump time, $\hat\tau\in ((k-1)h_n,kh_n)$, on a bin with $\Delta_k\widehat{[X,X]}>a^2 \vee u_n$, more precisely, we partition this bin into $R_n$ sub-intervals of lengths $(r_n+l_n)/n$ with $(r_n+l_n)$ an even integer. The jump window $(t_{l-l_n}^n,t_{l+r_n}^n)$, with length $h_n$ or smaller and $l=\lfloor \tau n\rfloor+1$, includes the price jump. 
$r_n$ determines the number of observations to the right of the price jump up to the end of the jump window, $l_n$ or $l_n+1$ is the number of observations to the left of the price jump down to the beginning of the jump window. For the price-jump estimation, we then cut out this jump window that contains $\tau$. Given $t^n_{l+r_n}$ and $t^n_{l-l_n}$, we thus use \eqref{bw} with the basis \eqref{spectral} centered around $Y_{t_{l+r_n}^n}-Y_{t_{l-l_n}^n}$ and with returns $\Delta_i^nY$ in a window $[t_{l-l_n}^n-h_n/2,t_{l-l_n}^n]$ to the left of the jump window and $[t_{l+r_n}^n,t_{l+r_n}^n+h_n/2]$ to the right of the jump window. This is the same as deleting observations $Y_{t_i^n}$ on $(t_{l-l_n}^n,t_{l+r_n}^n)$ and shifting observations $Y_{t_i^n}$ from the left and right towards the center. We identify a jump window by comparing $R_n$ pre-average jump estimators  
\begin{align}\label{argmax}\dot\iota=\operatorname{argmax}_{i=1,\ldots,R_n} \big|T^{LM}\big((k-1)h_n+(i-1/2)\frac{r_n+l_n}{n};\Delta_1^n Y,\ldots,\Delta_n^n Y\big)\big|\,,\end{align}
with the statistics from \eqref{lm}, averaging over $(r_n+l_n)/2\ll \sqrt{n}$ instead of $M_n$ observations. The final jump-size estimator is denoted by $\widehat{\Delta X}_{\hat\tau}$.
\newpage
\begin{prop}\label{argmaxprop}
When $R_n\rightarrow\infty$, with $R_n=\KLEINO(\sqrt{n})$, such that $(r_n+l_n)\propto n^{\delta}\rightarrow\infty$, for some $\delta>0$, the adjusted price-jump estimation using $\widehat{\Delta X}_{\hat\tau}$ with
\begin{align}t_{l-l_n}^n=(k-1)h_n+(\dot\iota-1)\frac{r_n+l_n}{n},~t_{l+r_n}^n=(k-1)h_n+\dot\iota\,\frac{r_n+l_n}{n}\,,
\end{align}
with $\dot\iota$ defined in \eqref{argmax}, satisfies \eqref{stablecltbwgen} in Proposition \ref{propbwgen}. 
\end{prop}
While $(r_n+l_n)$ is set by the econometrician, the two summands $l_n$ and $r_n$ are unknown and depend on the true value of $\tau$. Proposition \ref{argmaxprop} establishes asymptotically efficient price-jump estimation under noise, even if the jump time $\tau$ is unknown.

Analogously, the Lee-Mykland statistic \eqref{lm} can be adjusted for the unknown time point, $\tau$, in the jump window $(t_{l-l_n}^n,t_{l+r_n}^n)$. We estimate the price to the left of the jump window with $\hat P(t_{l+r_n}^n)$ and to the right of the jump window with $\hat P(t_{l-l_n-M_n}^n)$. This adjustment is also robust in the sense that Proposition \ref{proplm} remains valid.

Consider the two illustrative ``extreme'' examples in determining jump windows:
\begin{example}\onehalfspacing  $R_n=1$ implies cutting out the whole bin with $t_{l-l_n}^n=(k-1)h_n$ and $t_{l+r_n}^n=kh_n$. We can show that the price-jump estimator \eqref{bw} is consistent and preserves (almost) the optimal convergence rate in this case. However, the constant in the variance in \eqref{stablecltbw} increases when the jump window is of order $n^{-1/2}$. 
\end{example}
\begin{example}\onehalfspacing $R_n=nh_n-1$, when $\hat\tau_k=\operatorname{argmax}_i{\{t_i^n\in[(k-1)h_n,kh_n)|\,|\Delta_i^n Y|\}}$, implies centering \eqref{bw} around the largest absolute return on a bin. Since the noise is centered and its variance $\eta_{\tau_k}^2$ typically is rather small (see \cite{hans06}), the time of the largest absolute return might be considered a good candidate for the jump arrival and the method would require one fewer tuning parameter. In particular, if one addresses jumps much larger than $\eta_{\tau_k}$, the method could also perform well in practice. Theoretically, however, centering the jump window around the largest absolute return is only suitable if one assumes that $\eta_{\tau_k}\rightarrow 0$ when $n\rightarrow\infty$.  
\end{example}
\cite{aflwy16} point out that a central limit theorem for the DLE in the presence of market microstructure noise cannot generally be obtained with pre-averaging or related approaches. However, by focusing either on the \emph{tail} DLE, with some $a>0$ or assuming $r=0$ in Assumption \ref{jumps}, we derive the following asymptotic result:
\begin{prop}\label{propdle}\onehalfspacing Under Assumptions \ref{sigma}, \ref{jumps} and \ref{eta}, for any $a>0$ in that the L\'{e}vy measure $\mu$ does not have an atom, the estimator for the  DLE \eqref{lev} satisfies the feasible (self-scaling) central limit theorem
\begin{align}\label{cltdle}n^{\beta/2}\frac{\big(\widehat{[X,\sigma^2]}^d_1(a)-[X,\sigma^2]^d_1(a)\big)}{\big(\sum_{k=1}^{h_n^{-1}}\big(\widehat{\Delta X}_{\hat\tau_k}\big)^2\,8\hat\eta_{\hat\tau_k}\big(\hat\sigma^3_{\hat\tau_k}+\hat\sigma^3_{\hat\tau_k-}\big)\1_{\{\Delta_k\widehat{[X,X]}>a^2\vee u_n\}}\big)^{1/2}}\stackrel{(d)}{\longrightarrow} N(0,1)\,,\end{align}
with $\beta$ as in \eqref{beta}. If no price jump is detected, we set the estimate equal to zero. In particular, the limit theorem facilitates, for some $\alpha\in(0,1)$, an asymptotic level $\alpha$ test with asymptotic power 1 for testing the hypothesis $\tilde H_0:[X,\sigma^2]^d_1(a)=0$, against the alternative $\tilde H_1:[X,\sigma^2]^d_1(a)\ne 0$:
\begin{align}\varphi=\1_{\big\{|n^{\beta/2}\widehat{[X,\sigma^2]}^d_1(a)|>q_{1-\alpha/2}{\sqrt{\sum_{k=1}^{h_n^{-1}}(\widehat{\Delta X}_{\hat\tau_k})^2\,8\hat\eta_{\hat\tau_k}(\hat\sigma^3_{\hat\tau_k}+\hat\sigma^3_{\hat\tau_k-})\1_{\{\Delta_k\widehat{[X,X]}>a^2\vee u_n\}}}}\big\}}\,, \end{align}
where the $(1-\alpha/2)$ quantile of the standard normal law is denoted by $q_{1-\alpha/2}$.
\end{prop}
One loses no generality by imposing the scaling $T=1$; any fixed $T\in\mathbb{R}_+$ can be considered. The condition that the L\'{e}vy measure $\mu$ does not have an atom in $a$ is analogous to (10.76) in \cite{sahaliajacod}. There are only atoms in at most countably many values. According to \cite{sahaliajacod}, the condition holds for any $a>0$ as soon as $\mu$ has a density. This applies to all models used in finance with infinite jump activity. 
\begin{prop}\label{cordle}\onehalfspacing Under Assumptions \ref{sigma}, \ref{jumps} and \ref{eta} and under the specific case of finite jump activity, $r=0$ in Assumption \ref{jumps}, the estimator for the DLE,
\begin{align*}
\widehat{[X,\sigma^2]}^d_1=\sum_{k=2}^{h_n^{-1}-1}\widehat{\Delta X}_{\hat\tau_k}\big(\hat\sigma^2_{\hat\tau_k}-\hat\sigma^2_{\hat\tau_k-}\big)\,\1_{\{\Delta_k\widehat{[X,X]}> u_n\}}\,,\end{align*}
together with $\beta$ as in \eqref{beta} and $\varpi<\frac{1+\delta/2-1/4}{2+\delta/2}$ satisfies the feasible central limit theorem,
\begin{align}n^{\beta/2}\frac{\big(\widehat{[X,\sigma^2]}^d_1-[X,\sigma^2]^d_1\big)}{\big(\sum_{k=1}^{h_n^{-1}}\big(\widehat{\Delta X}_{\hat\tau_k}\big)^2\,8\hat\eta_{\hat\tau_k}\big(\hat\sigma^3_{\hat\tau_k}+\hat\sigma^3_{\hat\tau_k-}\big)\1_{\{\Delta_k\widehat{[X,X]}>u_n\}}\big)^{1/2}}\stackrel{(d)}{\longrightarrow} N(0,1)\,.\end{align}\end{prop}
The upper bound on $\varpi$ relates to Assumption \ref{eta} and the existence of higher moments of $\epsilon_t$. If all moments of the noise exist, the bound imposes no condition on the truncation. For $\delta\rightarrow 0$ in Assumption \ref{eta}, $\varpi<3/8$ leads to more conservative thresholds. Since $r=0$ in \eqref{beta}, we also derive the optimal rate  in this case. Although we  conjecture that this upper bound on $\varpi$ is not needed, it simplifies the proof considerably.

Proposition \ref{propdle} follows from combining our results on jump localization, the estimation of price jumps at detected jump times and from results of Corollary \ref{testvolajump} about volatility jump estimation. However, the proof cannot be extended in a similar way to the case $r\ne 0$ and $a=0$ when considering infinitely many small price jumps. It is unknown if an asymptotic distribution theory is possible in this general case. Propositions \ref{propdle} and \ref{cordle} give us exactly the statistics we require to apply in our data study, however.
\begin{remark}\onehalfspacing Propositions \ref{propdle} and \ref{cordle} indicate that, in the asymptotic results of the estimated DLE, the estimation error for the volatility jumps dominates the error for the price jumps. Consequently, the length of the jump window in Proposition \ref{argmaxprop} for price-jump estimation has asymptotically no effect on DLE estimation. Nevertheless, choosing $R_n>1$ is of interest from an applied point of view. Removing jump windows has a locally similar effect as downsampling the data to a lower observation frequency. Given the discussion by \cite{cop14} about spurious jump detection via downsampling, one would like to avoid deleting large jump windows in the empirical application. The refined method is superior to cutting out larger windows in that it poses less risk of estimating spuriously large jumps. 
\end{remark}
In addition to the DLE of \cite{aflwy16}, the leverage effect is also often defined in terms of a correlation statistic. To gain further insights across individual firms in the empirical Section \ref{sec:emp}, we follow \cite{jkm13} and consider a scaled measure of the DLE:  
 \begin{eqnarray} 
\frac{[X,\sigma^2]^d_T(a)}{\sqrt{[X,X]_T^d(a)[\sigma^2,\sigma^2]_T^d}(a)}
=\frac{\sum_{s\leq T}\Delta X_s\Delta\sigma_s^2\mathbbm{1}_{\{|\Delta X_s|>a \}}}{\sqrt{\sum_{s\leq T}(\Delta X_s)^2\mathbbm{1}_{\{|\Delta X_s|>a \}}}\sqrt{ \sum_{s\leq T}(\Delta \sigma^2_s)^2\mathbbm{1}_{\{|\Delta X_s|>a \}}}}\, , \label{cor}
\end{eqnarray} 
 that is, the correlation between contemporaneous price and volatility jumps. We may use $a=0$ in case of finite activity jumps, $r=0$ in Assumption \ref{jumps}. Note that \eqref{cor} is a path-wise defined, integrated measure. \eqref{cor} is a scalar parameter only under the restriction to time-homogeneous jump measures. Using Propositions \ref{propdle} and \ref{cordle}, and setting $T=1$, we obtain the following result:
\begin{cor}\label{corcorrelation}\onehalfspacing Under all conditions from Proposition \ref{propdle} and with
\begin{subequations}
\begin{align}\widehat{[\sigma^2,\sigma^2]}^d_1(a)&=\sum_{k=2}^{h_n^{-1}-1}\big(\hat\sigma^2_{\hat\tau_k}-\hat\sigma^2_{\hat\tau_k-}\big)^2\,\1_{\{\Delta_k\widehat{[X,X]}>a^2 \vee u_n\}}\,,\\ \widehat{[X,X]}^d_1(a)&=\sum_{k=2}^{h_n^{-1}-1}\big(\widehat{\Delta X}_{\hat\tau_k}\big)^2\,\1_{\{\Delta_k\widehat{[X,X]}>a^2 \vee u_n\}}\,,\end{align}
\end{subequations}
we derive a consistent estimator of \eqref{cor} with
\[\frac{\widehat{[X,\sigma^2]}^d_1(a)}{\sqrt{\widehat{[X,X]}^d_1(a)\widehat{[\sigma^2,\sigma^2]}^d_1(a)}}-\frac{[X,\sigma^2]^d_1(a)}{\sqrt{[X,X]_1^d(a)[\sigma^2,\sigma^2]_1^d}(a)}=\mathcal{O}_{\bar \P}(n^{-\beta/2})\,,\]  
with $\beta$ as in \eqref{beta}. Analogously, in the setup of Proposition \ref{cordle}, we obtain the same result for $a=0$.  
\end{cor}
\section{Simulations}\label{sec:sim}
This section reports the results of simulation studies of the finite-sample properties of the price-jump estimators, the corresponding price-jump tests and the discontinuous leverage statistics. The simulation study in \cite{bw16} evaluates the finite-sample inference on volatility jumps. 

This simulation study emulates that of \cite{lm12}. Although their theory only applies to the jump-diffusion setup, they simulate a more complex and realistic model, including stochastic volatility and time-varying noise. The efficient price follows
\begin{align}\label{xsim}X_t=1+\int_0^t\sigma_s\,dW_s\,,t\in[0,1],\end{align}
with Heston-type stochastic volatility,
\begin{align}\label{volasim}d\sigma_s^2=0.0162\,\big(0.8465-\sigma_s^2\big)\,ds+0.117\,\sigma_s\,dB_s\,,\end{align}
where $B$ and $W$ are two independent standard Brownian motions. We adopt the parameter values of \cite{lm12} in \eqref{volasim} and
assume 252 trading days per year and 6.5 trading hours a day. The model for the market microstructure noise is
\begin{align}\label{noisesim}\epsilon_{t_i^n}=0.0861\Delta_i^n X+0.06\big(\Delta_i^n X+\Delta_{i-1}^n X\big) U_i\,,\; i=0,\ldots,n,\end{align} 
with $(U_i)_{0\le i\le n}$ being a sequence of normally distributed random variables with mean 0 and variance $q^2$. We consider two parameterizations of $q$, which governs the noise level (market quality parameter). The cross-correlation between $X$ and noise violates one of our theoretical assumptions, but we expect no degradation in the performance of our approach. We estimate $q$ in the presence of serial correlation with the noise estimator suggested in Proposition 1 of \cite{lm12}.

We implement the self-scaling adaptive version of \eqref{bw} with pre-estimated optimal weights. The caption of Table \ref{tab1} gives values of $h_n$. The pre-averaging for the Lee-Mykland statistics \eqref{lm} refers to a block-size, $M_n=c\,\sqrt{n/k}$, where $k$ denotes the order of serial correlation in the simulated noise. The constant $c$ is chosen according to Table 5 of \cite{lm12}.\\[.2cm]
{\bf{\em Evaluation of the pre-average and spectral tests to infer price jumps}}\\[.2cm]
\cite{lm12} compare the performance of the noise-robust local jump tests in \cite{lm12} to those in \cite{lm8}, which are not designed to be robust to noise. We replicate this simulation study and compare the finite-sample performances of the statistics defined in \eqref{lm} and \eqref{bw}. Considering the power of the tests associated with Proposition \ref{proplm} (Lee-Mykland) and Proposition \ref{propbw} (our spectral method) allows us to compare our results to those in Table 4 of \cite{lm12}. Realizations of $Y_i=X_{t_i^n}+\epsilon_{t_i^n}$ are generated for one trading hour using time resolutions of 1, 2 and 3 seconds, respectively ($n=3600$, 1800, 1200). The jump size in $\tau$ is related to the noise level $q$, i.e., $\Delta X_{\tau}=0$ under the hypothesis and $\Delta X_{\tau}=q$, $2q$, $3q$ under the alternative.

\begin{table}[t]
\begin{center}
  \begin{threeparttable}
    \caption{Comparison of size and power of the two tests.\label{tab1}}
     \begin{tabular}{ccccccccc}
        \toprule
        \multicolumn{9}{c}{Moderate noise case, $q=0.0005$}\\
        \midrule
				Frequency ($n$) & \multicolumn{2}{c}{$\Delta X_{\tau}=0$}  & \multicolumn{2}{c}{ $\Delta X_{\tau}=q$} & \multicolumn{2}{c}{$\Delta X_{\tau}=2q$} & \multicolumn{2}{c}{$\Delta X_{\tau}=3q$}\\
				\midrule
				Test & LM & BNW & LM & BNW & LM & BNW & LM & BNW\\
				3 sec (1200)&0.049&0.045&0.199&0.274&0.473&0.677&0.777&0.924\\[-.2cm]
				&{\tiny{$(0.034)^*$}}&&{\tiny{$(0.059)^*$}}&&{\tiny{$(0.320)^*$}}&&{\tiny{$(0.786)^*$}}&\\
				2 sec (1800)&0.050&0.053&0.280&0.382&0.695&0.828&0.937&0.988\\[-.2cm]
				&{\tiny{$(0.030)^*$}}&&{\tiny{$(0.071)^*$}}&&{\tiny{$(0.483)^*$}}&&{\tiny{$(0.920)^*$}}&\\
				1 sec (3600)&0.049&0.056&0.281&0.594&0.697&0.982&0.950&1\\[-.2cm]
				&{\tiny{$(0.046)^*$}}&&{\tiny{$(0.091)^*$}}&&{\tiny{$(0.709)^*$}}&&{\tiny{$(0.988)^*$}}&\\
				\midrule
				\multicolumn{9}{c}{Large noise case, $ q=0.005$}\\
				\midrule
				Frequency ($n$) & \multicolumn{2}{c}{$\Delta X_{\tau}=0$}  & \multicolumn{2}{c}{ $\Delta X_{\tau}=q$} & \multicolumn{2}{c}{$\Delta X_{\tau}=2q$} & \multicolumn{2}{c}{$\Delta X_{\tau}=3q$}\\
				\midrule
				Test & LM & BNW & LM & BNW & LM & BNW & LM & BNW\\
				3 sec (1200)&0.052&0.049&0.296&0.996&0.803&1&0.997&1\\[-.2cm]
				&{\tiny{$(0.046)^*$}}&&{\tiny{$(0.275)^*$}}&&{\tiny{$(0.889)^*$}}&&{\tiny{$(0.997)^*$}}&\\
				2 sec (1800)&0.053&0.052&0.465&0.999&0.937&1&0.988&1\\[-.2cm]
				&{\tiny{$(0.046)^*$}}&&{\tiny{$(0.593)^*$}}&&{\tiny{$(0.998)^*$}}&&{\tiny{$(1)^*$}}&\\
				1 sec (3600)&0.050&0.049&0.829&1&0.994&1&0.997&1\\[-.2cm]
				&{\tiny{$(0.041)^*$}}&&{\tiny{$(0.918)^*$}}&&{\tiny{$(1)^*$}}&&{\tiny{$(1)^*$}}&\\
        \bottomrule
     \end{tabular}
    \begin{tablenotes}
      \small \item The table lists the simulated values of standardized test statistics \eqref{lm} and \eqref{bw}, from 6000 iterations for each configuration, exceeding the $0.05$-quantile of the standard normal. ``LM'' marks the Lee-Mykland test and ``BNW'' our proposed spectral test. We simulated from the model given by \eqref{xsim}, \eqref{volasim} and \eqref{noisesim}. In parentheses $(~)^*$, we report the values from Table 4 in \cite{lm12} of their analogous simulation study. According to Table 5 in \cite{lm12}, we used constants $c=1/19$ for $q=0.0005$ and $c=1/9$ for $q=0.005$ to determine $M_n$ in \eqref{lm} (for $\Delta X_{\tau}=2q, 3q$ and $q=0.005$, we doubled $M_n$, which increased the power). For \eqref{bw}, we used $h_n=\kappa\log{(n)}/\sqrt{n}$ with $\kappa \approx 5/12$ for $q=.0005$ and $\kappa\approx  2/3$ for $q=.005$.
    \end{tablenotes}
		\end{threeparttable}
	\end{center}
\end{table}
Table \ref{tab1} shows the simulation results, along with the values reported by \cite{lm12} in parentheses. Most of our results for the Lee-Mykland test closely track the ones reported by \cite{lm12}. Our results for the power under moderate noise and smaller jumps are a bit better than expected from \cite{lm12}, while some results in the large noise case are smaller. In the large noise case, we report values where $M_n$ is doubled compared to the constant adopted from \cite{lm12}, which led to higher power. The windows used for the spectral method are much larger than the values $M_n/n$ for the Lee-Mykland statistics. At first glance it might seem surprising that the power in Table \ref{tab1} increases for larger noise. This is not, however, because of large noise that makes precise testing and estimation more difficult, but because the jump sizes increase along with $q$. Large jumps naturally lead to better testing results. The size of both tests on the hypothesis appears to be accurate. The new spectral test \eqref{spectral} attains a considerably better power in all cases.\\[.2cm]
\newpage\noindent
{\bf{\em Evaluation of the pre-average and spectral estimators for price-jump sizes}}\\[.2cm]
In the same setup, we compare the performance of the jump-size estimators. Table \ref{tab2} confirms that, in all configurations, with the same optimal parameter choice as above, our spectral estimator attains a smaller root mean square error (RMSE) than the Lee-Mykland estimator. Efficiency gains are most relevant for the configuration with moderate noise and the smallest jump size. 
In this setup, our new estimator has a RMSE that is almost 50\% smaller for $n=3600$. For large noise and jump size $q$, our new estimator reduces the RMSE by 20\%. These significant improvements of estimation accuracy are particularly relevant because the moderate noise setting is realistic for current high-frequency data.

\begin{figure}
\centering
\includegraphics[width=5.7cm]{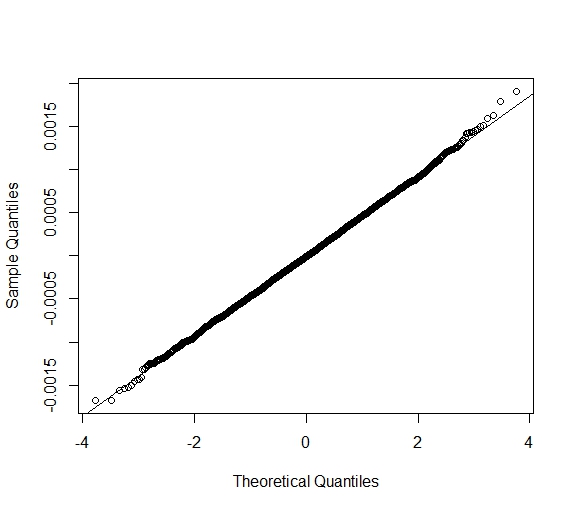}\hspace{-0.1cm}\includegraphics[width=5.7cm]{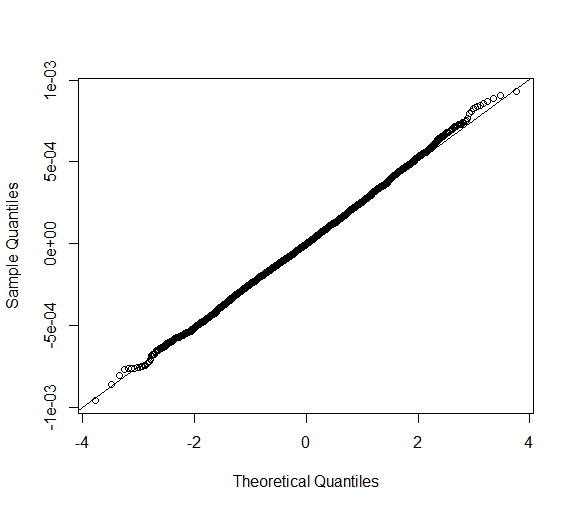} \\[-.5em]
\includegraphics[width=5.7cm]{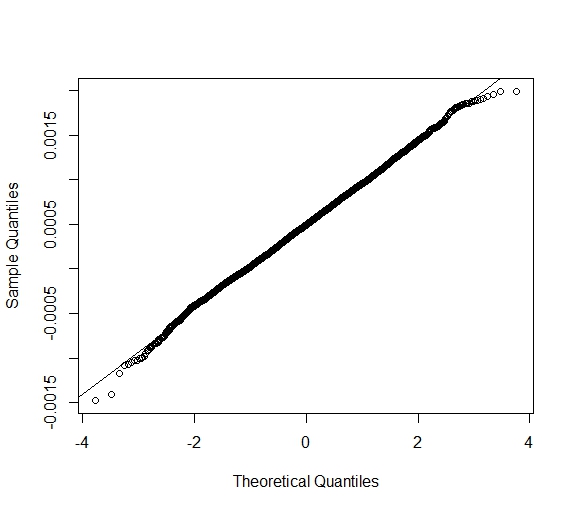}\hspace{-0.1cm}\includegraphics[width=5.7cm]{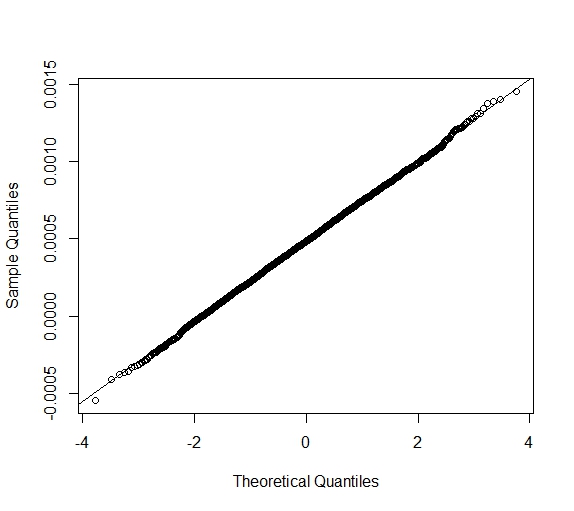} \\[-.5em]
 \caption{\label{qq}QQ-normal plots for the Lee-Mykland statistics (left) and our statistics \eqref{bw} (right). The top panels depict the 6000 iterations when $\Delta X_\tau=0$. The bottom panels show results for the iterations when $\Delta X_\tau=q=0.0005$.}
\end{figure}

Figure \ref{qq} demonstrates the finite-sample accuracy of the normal limit laws in \eqref{stablecltbw} and \eqref{stablecltlm}. The empirical distributions closely approximate their normal asymptotic limit.\\[.2cm] 
\begin{table}[t]
\begin{center}
  \begin{threeparttable}
	\caption{Comparison of RMSEs for the two price-jump size estimators.\label{tab2}}
     \begin{tabular}{ccccccc}
        \toprule
        \multicolumn{7}{c}{Moderate noise case, $q=0.0005$}\\
        \midrule
				Frequency ($n$) ~~~& \multicolumn{2}{c}{ $\Delta X_{\tau}=q$} & \multicolumn{2}{c}{$\Delta X_{\tau}=2q$} & \multicolumn{2}{c}{$\Delta X_{\tau}=3q$}\\
				\midrule
				Estimator  & LM ~~& BNW~~~ & LM ~~& BNW~~~ & LM ~~& BNW~~~\\
				3 sec (1200)&11.0&9.9&11.1&10.2&11.9&10.8\\
				2 sec (1800)&6.8&5.3&6.9&6.0&7.9&6.8\\
				1 sec (3600)&4.7&2.6&4.8&3.6&6.3&4.7\\
				\midrule
				\multicolumn{7}{c}{Large noise case, $q=0.005$}\\
				\midrule
				Frequency ($n$) & \multicolumn{2}{c}{ $\Delta X_{\tau}=q$} & \multicolumn{2}{c}{$\Delta X_{\tau}=2q$} & \multicolumn{2}{c}{$\Delta X_{\tau}=3q$}\\
				\midrule
				Estimator & LM & BNW & LM & BNW & LM & BNW\\
				3 sec (1200)&14.8&14.4&15.0&14.5&15.2&14.5\\
				2 sec (1800)&10.0&9.4&10.2&9.5&10.6&9.5\\
				1 sec (3600)&5.6&4.5&5.9&4.6&6.4&4.6\\
        \bottomrule
     \end{tabular}
			\begin{tablenotes}
      \small \item The table lists the root mean square errors, rescaled by $10^4$, of the estimators \eqref{lm} and \eqref{bw}, from 6000 iterations for each configuration under the alternative when price jumps are present. ``LM'' marks the Lee-Mykland estimator and ``BNW'' our proposed estimator. We simulated from the model given by \eqref{xsim}, \eqref{volasim} and \eqref{noisesim}. Tuning parameters are reported in Table \ref{tab1}.
  \end{tablenotes}
	  \end{threeparttable}
		\end{center}
\end{table}
{\bf{\em Evaluation of the discontinuous leverage estimator}}\\[.2cm]
\begin{figure}[ht]
\centering
\includegraphics[width=5.7cm]{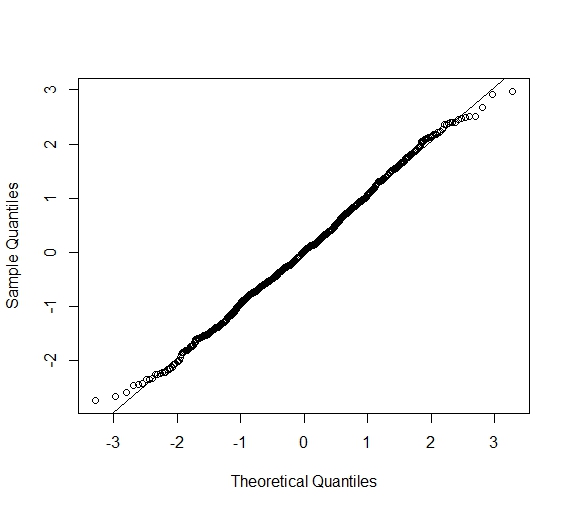}\hspace{-0.1cm}\includegraphics[width=5.55cm]{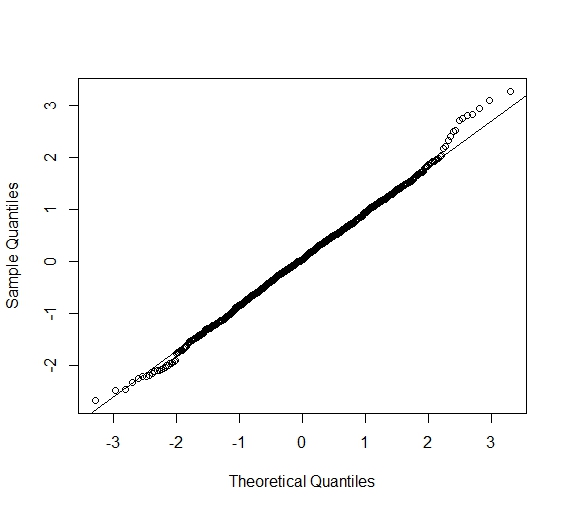} 
 \caption{\label{qq2}QQ-normal plots for \eqref{cltdle}, 1-second (left) and 3-second frequencies (right).}
\end{figure}
We modify the simulation setup by adding one jump at a random time to the volatility in \eqref{volasim}. The volatility jump size is set to the median value from the empirical sample described in Table \ref{desc1}. To create discontinuous leverage, we implement a contemporaneous downward price jump of 0.2\%, which is comparable to the sizes in Figure \ref{Fig2}. Using a rather large price jump and average volatility jump-size allows us to study the finite-sample accuracy of the result  \eqref{cltdle}. We can analyze the DLE estimator \eqref{lev} because thresholding reliably detects such jumps. We simulate one trading day with observation frequencies of 1, 2 and 3 seconds, frequencies that generate 23400, 11700 and 7800 observations, respectively, over the day. 

We estimate the DLE in a model with moderate microstructure noise. We first estimate spectral statistics over a partition of the whole day, identifying price jumps by thresholding. Next, we estimate the squared volatility before and after the jump by local averages of the bin-wise, parametric estimates over 8 bins. Then, we estimate the local jump size using \eqref{bw} and implement the refinement from Section \ref{sec:3.3} for unknown jump times. We partition the bin with the detected jump in $R=6$ equidistant sub-intervals and apply the adjusted jump size estimation using \eqref{argmax}. The window sizes for the first step and the price-jump estimation are equal: we use $h^{-1}=100$ for 1-second frequency and $h^{-1}=50$ for the two smaller frequencies. The spectral cut-off frequency is set to $J=30$ in all cases. Estimates are reasonably robust to different values of $h$ and $J$.

For the fixed, true value -$2.324$ of the DLE \eqref{dle},\footnote{We rescale all DLE values with $10^7$.} we obtain these estimates:\\
\hspace*{-.2cm}\begin{tabular}{lccc}
Frequency& 1 sec & 2 sec & 3 sec\\
Bias&-0.04 & -0.02 & -0.03\\
Variance & 0.16 & 0.19 & 0.21\,
\end{tabular}\\[-.1cm] 

The inherently slow convergence rate of the estimation leads to pronounced finite sample variances. 
Figure \ref{qq2} shows QQ-normal plots for the test statistic obeying the central limit theorem \eqref{cltdle}. The normal distribution fits reasonably well for all frequencies. Our test for the DLE attains very high power (approximately 99\%) in the case of one observation per second and only slightly smaller power for the lower observation frequencies. Overall, simulations indicate that the estimation performs well in this complex environment.
\clearpage 
\section{The discontinuous leverage effect in stock prices}\label{sec:emp}
This section presents results of applying the spectral methods of Section \ref{sec:3} to stock price data. We first introduce the dataset and discuss how to estimate price and volatility jumps on these data.  Second, we investigate the DLE, i.e., a covariation measure, and the correlations of price and volatility cojumps. Finally, we explain the cross-sectional variation in the DLE and correlation estimates across firms. 
\subsection{Price and volatility cojumps}
We use NASDAQ order book data from the LOBSTER database. Initially, we pick the 30 stocks with the largest  market capitalizations from each of the 12 NASDAQ industries for a total of $12*30=360$ stocks.\footnote{The industries can be found on www.nasdaq.com/screening/industries.aspx. The year 2013 serves as the baseline year.} The sample spans January 1, 2010 to December 31, 2015, 1,509 days with trading from 9:30 to 16:00 EST. The tick-by-tick data shows evidence of market microstructure noise in that, for instance, returns have significant negative first-order autocorrelation. The test of \cite{ax16}, equation (40), displays significant noise for 50\% of all stocks, across all trading days.\footnote{\cite{ax16} report similar percentages for the S\&P100 in their Table 4. To control the overall significance level of tests across firms and trading days, we use the \cite{bh95} step-up procedure at level $\alpha=0.1$. In case of no market microstructure noise our methods remain valid.} As shown in the simulations of \cite{wbl16}, spectral estimators perform particularly well with liquid stocks, i.e., those having at least about one trade every 15 seconds. To restrict the analysis to very liquid stocks, we exclude trading days for a given stock with fewer than 1,500 trades.\footnote{Results are robust to  higher (2,000) and lower (1,000) thresholds.} This selection procedure reduces the number of firms to 320. We focus on transactions with non-zero returns but do not adjust the data further; that is, we do not clean or synchronize trades. The number of observed trades varies substantially across stocks and days. There is a maximum of 227,139 intradaily observations for the Apple Inc.\ stock on September 9, 2014; the median number of daily transactions across stocks is much smaller, only 5,977. 

The local jump detection and estimation takes the time-varying trading activity into account. We partition each trading day $d$ into $h^{-1,(d,s)}=K^{(d,s)}$ bins for every stock $s$. As suggested by our theoretical results, the number of bins $k=1,...,K^{(d,s)}$ grows with the number of trades $n^{(d,s)}$ with $K^{(d,s)}=\lfloor 3\sqrt{n^{(d,s)}}\log(n^{(d,s)})^{-1}\rfloor$. We detect price jumps by applying the adaptive threshold, $ u_k^{(d,s)}=2\log(K^{(d,s)})/K^{(d,s)}\hat\sigma^{2,(d,s)}_{k,pil}$, to bin-wise quadratic variation estimates \eqref{quadad}. It is well-known that the number of detected price jumps depends on the thresholding procedure in the sense that a lower threshold usually increases the number of detected small price jumps.\footnote{The main results about the DLE are robust against different threshold choices. As a robustness check, we substitute the log($K$) term of the threshold to log(log($K$)), which increases the number of price-jump days per stock from around 14\% to 29\%.} We find that relatively small volatility changes at price-jump time points strongly influence the DLE estimates. For that reason, we apply the test for volatility jumps, as proposed by \cite{bw16}, to focus on price jumps with significant contemporaneous volatility jumps.\footnote{Note that in \eqref{lev}, summands without volatility jumps ``automatically'' cancel out because $\sigma_{\tau}^2=\sigma_{\tau-}^2$. To control the overall significance level of tests across firms and price jumps at level $\alpha=0.1$,  we use the  \cite{bh95} step-up procedure (the false discovery rate).}
The tests for volatility jumps reduce the influence of the price-jump-detection threshold on the DLE estimates.
For price-jump estimation, we partition jump bins into $R=6$ sub-intervals and center \eqref{bw} around the cut-out return obtained via \eqref{argmax}.\footnote{Note that centering the jump estimator around the largest absolute return on a detected bin, as described in Example 2, does not change the main conclusions about the DLE. However, individual estimates of price-jump sizes can differ quite substantially.} The number of frequencies studied on each bin is $J^{(d,s)}=5\log(n^{(d,s)})$. We average the truncated spectral statistics over $R^{(d,s)}=\lceil3 \sqrt[4]{n^{(d,s)}}/\log(n^{(d,s)})\rceil$ bins to estimate spot volatility to the right and left of the detected price jump.  
\begin{figure}[t]
\centering
\includegraphics[width=13cm]{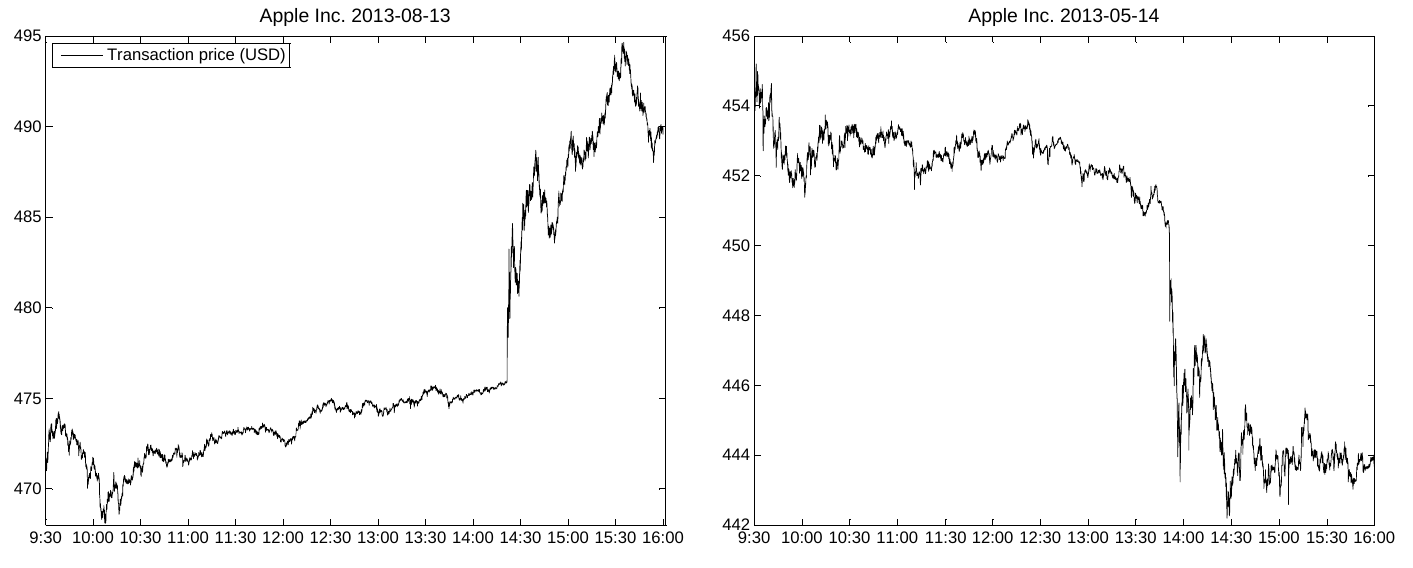}
\caption{\label{Fig2} Price process at the NASDAQ stock exchange of Apple Inc.\ on two different days with price-volatility cojumps. Number of trades: 87,445 (left), 40,707 (right). }
\end{figure} 

Figure \ref{Fig2} shows two examples of price-volatility cojumps of the Apple Inc.\ stock, an upward price jump in the left panel and a downward price jump in the right panel.  The estimates of the price jumps are $ 0.27\%$ and -$0.24\%$, respectively. Note that if one would approximate the price-jump sizes just by looking at Figure \ref{Fig2} and assuming a small noise level, one may expect much larger price-jump estimates. 
\cite{cop14} and \cite{bhls09} explain that  seemingly large returns often consist of smaller, unidirectional returns on a short time interval.\footnote{While \cite{cop14} attribute a local drift to such phenomena, \cite{bhls09} explain this characteristic by the microstructure of the orderbook and call it ``gradual jumps''.}   
This explains how downsampling to lower observation frequencies can affect both jump detection and the estimation of price-jump sizes. 
Figure \ref{Fig2} also suggests that volatility jumped contemporaneously with the price jump. That is, the variability of the stock price appears in both cases much smaller before the price jump than afterwards. This apparent jump in volatility is not directly determined by the price jump mechanically feeding through to higher volatility. Indeed estimated changes in volatility only use log-price information from bins that neighbor the price-jump bin. The increase in spot volatility evaluated approximately 30 minutes before and after the price jumps is $184\%$ (left panel of Figure \ref{Fig2}) and $163\%$ (right panel of Figure \ref{Fig2}). Note that the strong upward jumps in both the price and volatility processes, in the left panel of Figure \ref{Fig2}, is not consistent with the negative price-volatility cojump correlation in high-frequency data that \cite{br16} report for S\&P 500 futures.
\begin{table}
\centering \small
\caption{Price and volatility cojumps: NASDAQ order book, 2010-2015.}
\label{desc1}
\begin{tabular}{lcccccccccccc}  \hline\hline \vspace{-.8em} \\[-.25em]
\multirow{2}{0.5cm}{Conditioning criteria}&&\multirow{2}{1.8cm}{\centering \# of cojumps}  &\multicolumn{3}{c}{Price-jump size}&&\multicolumn{3}{c}{Volatility-jump size}  \\  \cline{4-6}\cline{8-10} \\[-1.05em]
&&&$Q_{0.25}$&$Q_{0.5}$&$Q_{0.75}$&&$Q_{0.25}$&$Q_{0.5}$&$Q_{0.75}$ \\ \hline \\[-.6em]
\multicolumn{6}{l}{{\bf Panel A:  Apple Inc.\ stock}} \\ \hline
All jumps& & 209& -0.095&-0.033&0.048&&43.8&88.6&193.4\\
\multicolumn{4}{l}{Positive price jumps}&		\\
$\cdot$ All&&83& 0.029&0.070&0.152&&-28.4&84.8&180.3\\
$\cdot$ Market  &&20&0.073&0.152&0.266&& -28.2&124.7&266.0\\
$\cdot$ Idiosyncratic  &&63&0.028&0.061&0.119&& -28.4&80.6&153.3 \\
\multicolumn{4}{l}{Negative price jumps}&		\\
$\cdot$ All&&  126& -0.137&-0.084&-0.046&& 51.4&89.6&216.9 \\
$\cdot$ Market&& 19& -0.172&-0.110&-0.050&& 86.2&316.5&504.9 \\
$\cdot$ Idiosyncratic &&107&-0.130&-0.084&-0.046&&49.4&82.5&174.4 \\[.35em] 
\multicolumn{6}{l}{{\bf Panel B:  Mean across all stocks}} \\ \hline
All jumps&& 73.8& -0.115&0.014&0.152&&39.2&137.3&299.1\\
\multicolumn{4}{l}{Positive price jumps}&		\\
$\cdot$ All& &38.5& 0.108&0.175&0.286&&-23.3&114.4&290.7\\
$\cdot$ Market&  &9.9& 0.139&0.206&0.327&&154.2&269.5&361.9 \\
$\cdot$ Idiosyncratic& &28.6& 0.101&0.164&0.266&&-21.8&108.9&230.1 \\
\multicolumn{4}{l}{Negative price jumps}&		\\
$\cdot$ All& &35.4& -0.254&-0.160&-0.101&&76.7&145.6&315.0 \\
$\cdot$ Market& &7.6&-0.283&-0.194&-0.141&& 168.8&276.8&444.6  \\
$\cdot$ Idiosyncratic& &27.9&-0.242&-0.154&-0.098&&74.7&119.5&304.2 \\
\hline \hline  
\multicolumn{10}{p{12.45cm}}{\footnotesize Notes: Quantiles ($Q$) of the jump distributions are in percent. Market jumps refer to days with jumps in the NASDAQ composite index. Idiosyncratic jumps refer to days without jumps in the NASDAQ composite index.}
\end{tabular}
\end{table}

To get deeper insights about price and volatility cojumps, Table \ref{desc1} shows summary statistics for detected cojumps and quantiles of the respective jump distributions. Panel A of Table \ref{desc1} shows summary statistics for the Apple Inc.\ stock; Panel B displays averages across the 320 stocks. We condition results on the sign of price jumps and whether they are market jumps or idiosyncratic. Following \cite{ltt17}, we use a market index to observe systematic risk. The market proxy is the NASDAQ Composite Index, which is the market capitalization-weighted index of about 3,000 equities listed on the NASDAQ stock exchange. The price-jump detection method proposed in Section \ref{sec:3.3} finds market jumps on 8\% of the days in the sample. We define idiosyncratic  jumps as discontinuities where the index displays no contemporaneous market jump. 

The top row of Panel A of Table \ref{desc1} shows that the Apple Inc.\ stock price displays 209 contemporaneous price-volatility cojumps, with more downward price jumps (126 or 60\%) than upward price jumps (83 or 40\%). Panel B of Table \ref{desc1} shows that the average number of price-volatility cojumps in the six-year sample across all individual stocks is 73.8. Columns three to five and six to eight of Table 3 show the quantiles of the price-jump and volatility-jump distributions. They indicate that idiosyncratic jumps are smaller than market jumps. The  magnitude of price jumps is in line with the sizes of -$0.15$ to 0.18\% reported by \cite{lm12} for the IBM stock in 2007. 

The magnitude of volatility jumps is striking. The 0.75 empirical quantile of the volatility-jump distribution of the Apple Inc.\ stock for negative market price jumps is about 505\%. That is, volatility frequently jumps to more than five times its pre-jump size when prices jump down. The analogous 0.75 quantile for volatility jumps, conditional on a  negative market price jump, averaged across all firms is 445\%. Scheduled news announcements are known to reduce trading and volatility right before the announcement but portend a strong response afterwards, which is manifested in large volatility jumps. 
The rows labeled ``market,'' in Panel B of Table \ref{desc1}, show that volatility jumps are usually positive for both positive and negative price jumps. Overall, the volatility-jump distribution is right-skewed, indicating the important role of upward jumps in volatility.  

\subsection{The discontinuous leverage effect}
\begin{figure}[t]
\centering
\includegraphics[width=11cm]{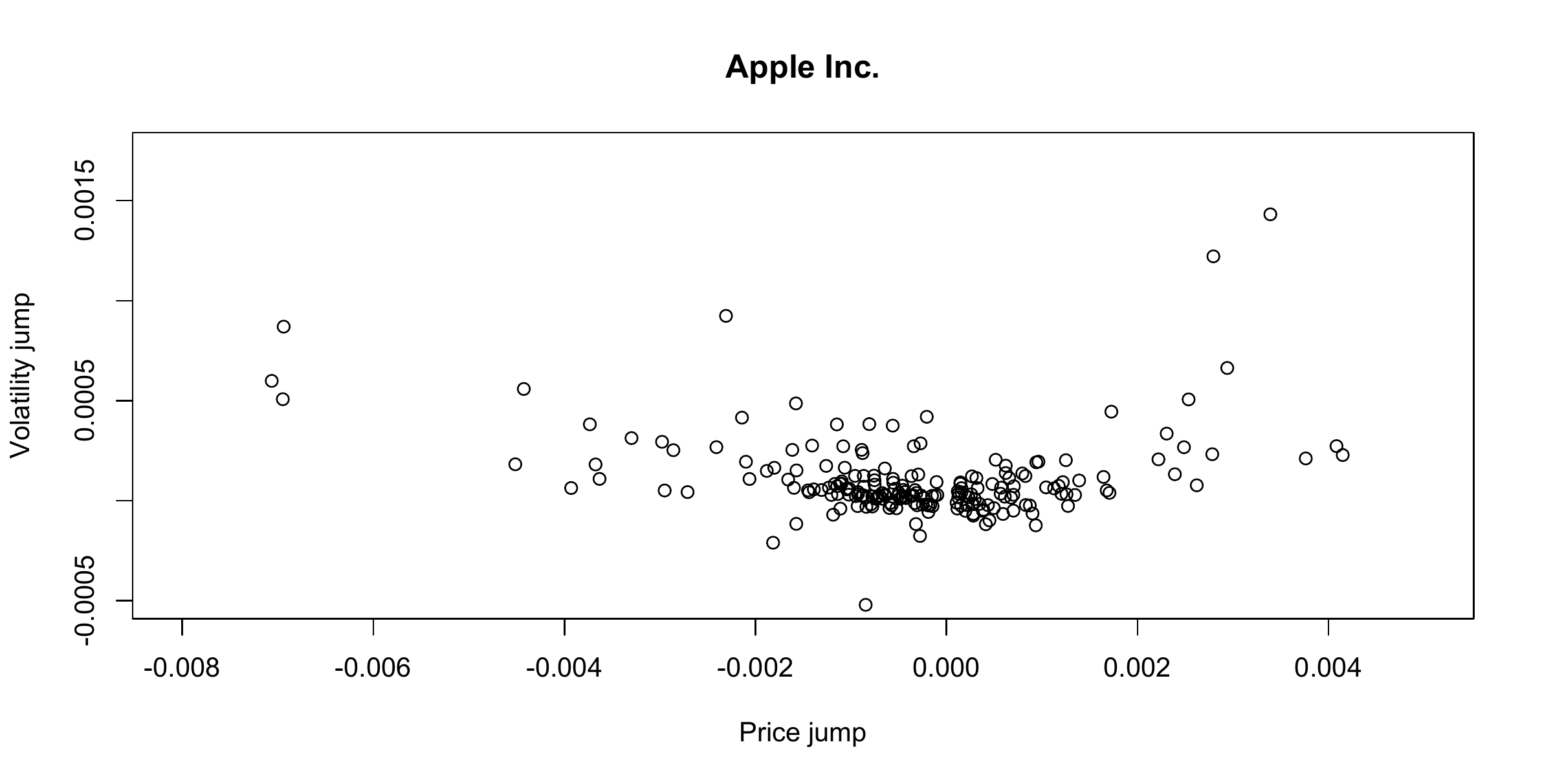}\hspace{-0.1cm}
\caption{Scatter plot of price jumps and contemporaneous volatility jumps. Sample period 2010-2015.} 
\label{scatter}
\end{figure}
This subsection characterizes the DLE of contemporaneous price and volatility jumps. 
Figure \ref{scatter} illustrates a typical relation between price jumps and contemporaneous volatility jumps using data from Apple Inc.\  from 2010 to 2015. Following \cite{dps00} and \cite{br16}, one would expect an unconditional, negative linear relation between the price and volatility-jump sizes. However, the figure does not depict such a uniformly negative relation. Row 1 of Table \ref{results} 
documents the absence of an unconditional price-volatility cojump relation across firms. That is, the  test \eqref{cltdle} rejects the null hypothesis of no DLE for only 10\% of the 320 firms. The DLE estimates and correlation \eqref{cor} are usually close to zero, with inconsistent signs across firms. The median DLE across all firms is $0.17$; the corresponding correlation is 0.01. In other words, there is no prevalent unconditional leverage effect using either measure of leverage. 

This result confirms previous negative findings of parametric asset pricing models by \cite{cggt03}, \cite{ejp03} and \cite{e04}, who use U.S.\ stock index and option data. \cite{jkm13} also find no significant correlation of price-volatility cojumps in one-minute S\&P 500 ETF data. Row 1 of Table \ref{results} thus extends the literature's negative results on discontinuous leverage to the cross-section of individual stock price processes.
\begin{table}
\begin{center} \small
\caption{The discontinuous leverage across NASDAQ firms. }
\label{results}
\begin{tabular}{clccccccccccc}  \hline\hline \vspace{-.8em} \\ 
\multirow{2}{*}{Row}&Conditioning&Rejection&&\multicolumn{3}{c}{DLE quantiles}&&\multicolumn{3}{c}{Correlation quantiles } \\ \cline{5-7} \cline{9-11} 
&criteria&rate&& $Q_{0.25}$&$Q_{0.50}$&$Q_{0.75}$&& $Q_{0.25}$&$Q_{0.50}$&$Q_{0.75}$\\
\hline \vspace{-.8em} \\ 
1&All jumps&  0.10&&-0.37&0.17&1.01&&-0.12&0.01&0.15\\[.3em]
&Positive price jumps&		\\
2&$\cdot$ All&  0.65&&1.14&2.08&3.87&		& 0.11&0.28&0.44\\
3&$\cdot$ Market  & 0.89&&1.64&2.88&4.48&				& 0.12&0.32&0.55 \\
4&$\cdot$ Idiosyncratic  & 0.66&& 0.84&1.63&3.15&						&0.09&0.26&0.47\\[.3em]
&Negative price jumps&		\\
5&$\cdot$ All&  0.63&&-3.72&-1.75&-0.88&	& -0.47&-0.29&-0.09\\
6&$\cdot$ Market& 0.85&&-3.75&-2.02&-1.03&& -0.62&-0.31&-0.15\\
7&$\cdot$ Idiosyncratic & 0.67&&-3.53&-1.69&-0.92&& -0.52&-0.32&-0.14\\[.25em] \hline
\hline 
\multicolumn{11}{p{13.6cm}}{\footnotesize Notes: The rejection rate indicates the percentage of firms having a significant DLE. We control the overall significance at level $\alpha=0.1$ with the step-up procedure of \cite{bh95}. DLE quantiles refer to a firm's average DLE, rescaled by $\times10^{7}$. The empirical quantiles contain all DLE estimates and correlation coefficients $Corr(\Delta X,\Delta\sigma^2)$ across firms.  }
\end{tabular}
\end{center}
\end{table}

Given that it is difficult to reject the hypothesis of no DLE, the question arises if we should expect the discontinuous relation to be similar to that of the continuous leverage.
As discovered by \cite{lln11}, specific events cause large jumps and those jumps are relatively rare. Volatility jumps are very large on impact, but the level of volatility often subsequently decays quickly in the direction of a pre-event level. The impact of news potentially drives common price and volatility jumps, as described by 
\cite{pv12}. We conjecture that such news effects usually trigger upward jumps in volatility, regardless of the effect on prices, and thus produce a positive (negative) correlation of volatility jumps with contemporaneous upward (downward) price jumps. To investigate this response pattern, we condition the DLE estimates on the signs of the price jumps. 
 
Rows 2 and 5 of Table \ref{results} show the outcomes of the DLE test \eqref{cltdle} conditional on upward and downward price jumps, respectively. We focus on stocks with more than 10 price-volatility cojumps and exclude jumps larger than six standard deviations, which leaves us with 307 firms. Quantiles of the DLE estimates and correlations indicate that the DLE is negative for downward price jumps and positive for upward price jumps. That is, the leverage statistic quantiles are uniformly positive (negative) for positive (negative) price jumps. 
Row 2 of Table \ref{results} shows that 65\% of the firms display a significant DLE if prices jump up. Similarly, row 5 of Table \ref{results} shows that 63\% of the firms have a statistically significant DLE for negative price jumps. The positive (negative) relation between positive (negative) price jumps and contemporaneous volatility jumps is also visible in the scatter plot in Figure \ref{scatter}. 

In addition to conditioning on the sign of the price jump, we consider the fact that standard asset pricing models price different sources of risk differently. Systematic jumps are often  related to macroeconomic news announcements and trigger cojumps across a large fraction of all firms while firm-specific jumps likely reflect idiosyncratic risk. 

Conditioning on whether price jumps are market-wide or idiosyncratic reveals a strong conditional relation between discontinuities in prices and volatility (see rows 3 and 6 of Table \ref{results}). We focus on firms having more than 10 market price-volatility cojumps and omit jumps larger than six times its standard deviation. This shrinks the number of firms to 230. For this sample,  market upward jumps and contemporaneous volatility jumps (see row 3 of Table \ref{results}) display a significant DLE for 89\% of the firms. The median DLE estimate across all firms for a single price-volatility cojump is $2.88$. The median correlation between positive market jumps and volatility jumps is 0.32. Downward market jumps (see row 6 of Table \ref{results}) exhibit a significant downward sloping relation for 85\% of the firms. The median DLE estimate for a single price-volatility cojump is -$2.02$ with a corresponding correlation of -$0.31$. 
A comparison of rows 3 to 4 and 6 to 7 of Table \ref{results} shows that market jumps are usually more strongly correlated with contemporaneous volatility jumps than are idiosyncratic price jumps. Market jumps show a stronger conditional DLE than do idiosyncratic jumps because market events coincide with large price and volatility cojumps. This allows us to conclude that the tail DLE is particularly strong. 

In contrast to market jumps, idiosyncratic jumps are smaller, coming more from the center of the jump distributions, and display a weaker DLE. Rows 4 and 7 of Table \ref{results} indicate that about 66\% of the stocks have a significant DLE for idiosyncratic jumps. 

All in all, two forces prevent an unconditionally negative DLE: First, the sign of the price-volatility cojump depends on the sign of the price jump. That is, positive (negative) price jumps are positively (negatively) correlated with contemporaneous volatility jumps. Second, the DLE is stronger for market price jumps and more often nonsignificant for idiosyncratic price jumps.

The positive (negative) correlation between upward (downward) price jumps and contemporaneous volatility jumps might explain why \cite{jkm13} find no significant, unconditional correlation between price and volatility jumps, while \cite{tt11} report a strong positive relation between squared price jumps and jumps in volatility. Our results indicate that one would expect a positive unconditional DLE between squared price jumps and volatility. The weaker relation between idiosyncratic price jumps and volatility jumps relates to \cite{y12}, who model a time-varying leverage effect in a (semi)parametric stochastic volatility model where the time-variation is associated with the size of returns. By conditioning on positive and negative price jumps, we focus on Yu's positive and negative extreme states. Our analysis indicates that it is important to distinguish market jumps and idiosyncratic jumps, which roughly implies distinguishing the tail from the rest of the price-jump distribution.
\subsection{Drivers of the discontinuous leverage}  
This subsection investigates several variables that might drive the cross-sectional variation in the DLE. A prime candidate for such a causal factor is a firm's debt-to-equity ratio. This explanatory variable is due to the original interpretation of the continuous leverage effect as stemming from the levered nature of equity.  
 The original interpretation of the leverage effect is that a lower stock price will increase the debt-to-equity ratio and make the stock price more volatile. Therefore, we investigate a firm's debt-to-equity ratio as a cross-sectional driver of the DLE. 

We regress the DLE and correlation estimates for market upward and downward price jumps on firms' debt-to-equity ratios and the volatility levels.\footnote{The volatility level can be considered as a component of financial leverage; see \cite{br12} for example.} 
We treat the positive and negative DLEs and correlations in separate regressions as an empirical matter because a given firm usually has asymmetric positive and negative DLEs.
Following the asset pricing literature, we include common firm characteristics, such as firm size, the book-to-market ratio, the book value of equity, profit, the cash-income ratio, the price-earnings ratio, and a liquidity measure (number of trades) to control for other influences. Compustat provides firm characteristics; we average characteristics and trades over the 2010-2015 sample before using them in the cross-sectional regression. 
\begin{table}
\begin{center} \small
\caption{Regressions with discontinuous leverage estimates and correlations.  } \label{reg}
\begin{tabular}{lccccccccc}  \hline\hline \vspace{-.8em} \\ 
Dependent variable& \multicolumn{4}{c}{Explanatory variables}&$ R^2$\\ \cline{1-1} \cline{2-7} \\[-.7em]
&Debt-to-equity&Vola level&Size &\# Trades&   \\ \cline{2-5} \\[-1.6em]
\multicolumn{4}{l}{Leverage $\widehat{[X,\sigma^2]}^d_T$} \\[.2em]
$\cdot$ Positive price jumps&$\underset{(0.12)}{-0.523^\dagger}$ &$\underset{(0.00)}{\bf 0.081}$&$\underset{(0.13)}{2.910^\ddagger}$&$\underset{(0.31)}{ -0.298^\ddagger}$&0.13\\
$\cdot$ Negative price jumps&$\underset{(0.50)}{0.293^\dagger}$ &$\underset{(0.03)}{\bf -0.015}$&$\underset{(0.24)}{ -2.950^\ddagger}$&$\underset{(0.05)}{\bf -0.385^\ddagger}$&0.12\\ 
\\[-1.1em]
\multicolumn{4}{l}{Correlation $Corr(\Delta X,\Delta \sigma^2)$}  \\[.2em] 
$\cdot$ Positive price jumps& $\underset{(0.33)}{-0.100}$ &$\underset{(0.00)}{\bf - 876.1}$&$\underset{(0.70)}{ -2.890^\dagger}$&$\underset{(0.96)}{-1.220^\dagger}$&0.08\\ 
$\cdot$ Negative price jumps& $\underset{(0.04)}{\bf -0.401}$ &$\underset{(0.02)}{\bf 885.9}$&$\underset{(0.64)}{4.300^\dagger}$&$\underset{(0.83)}{{-1.500}^\dagger}$&0.07\\[.2em]


\hline \hline 
\multicolumn{9}{p{11.8cm}}{\footnotesize Notes: Regressions refer to the market jumps. $\dagger:=\times 10^{-7}$, $\ddagger:=\times 10^{-11}$.  
$p$-values in parentheses. Cross-sectional regressions include further firm characteristics (not shown). Sample of 215 firms.}
\label{regre}
\end{tabular}
\end{center}
\end{table}

The first two rows in Table \ref{reg} show the debt-to-equity ratio does not significantly influence the cross-sectional variation in the DLE estimates. However, row 4 of Table \ref{reg} shows that firms with a higher debt-to-equity ratio have a significantly stronger correlation of negative market-wide price-volatility cojumps. An increase in the debt-to-equity ratio of 10 percentage points decreases the negative correlation of price and volatility cojumps by 0.04. This is consistent with the traditional interpretation of the leverage effect.

In contrast to the debt-to-equity ratio, the regressions show that the level of integrated volatility consistently explains the cross-sectional variation in DLE estimates across specifications. Firms with higher integrated volatility have a larger DLE (in absolute terms). This finding is consistent with the idea that systematic risk is not diversifiable. It implies that firms that respond more strongly to systematic risk, i.e.\ market jumps, tend to have higher volatility levels, hence are more risky. In contrast to the DLE estimates, the price-volatility cojump correlations decrease in absolute terms with the level of volatility.\footnote{The sign of the coefficient on the volatility regressor depends on whether the regressand is the covariance or the correlation. This sign reversal may occur because the correlation is the DLE scaled by the product of price-jump and volatility-jump variation and particularly the volatility-jump variation is highly positively correlated with the volatility level.} A relatively high level of volatility results in a more dispersed relation between price jumps and contemporaneous volatility jumps. This result seems to contradict the idea that the debt-to-equity ratio should have a positive relation with the volatility level, see \cite{br12} for example, and once more indicates the different natures of the continuous and discontinuous leverage.

\section{Conclusion}
This paper makes both methodological and empirical contributions to the literature on contemporaneous price and volatility jumps. 
We propose a nonparametric estimator of the discontinuous leverage effect (DLE) in high-frequency data that is robust to the presence of market microstructure noise. The new estimator allows us to study transactions data from the order book without down-sampling to a lower, regular observation frequency. For DLE estimation, we develop an efficient jump estimator for unknown jump times. We document the estimator's superior asymptotic and finite sample qualities compared to a method with pre-average jump-size estimation.  

Previous research has found it difficult to empirically document a DLE. Studying contemporaneous price and volatility jumps of 320 individual NASDAQ stocks from 2010 to 2015, we also find mixed and mostly insignificant, unconditional DLEs when considering all detected price and volatility cojumps. We show that the event-specific nature and distinct sources of jumps obscure the true relation between price and contemporaneous volatility jumps. We establish that a strong and significant DLE exists by conditioning on the sign of price jumps and on whether the price jumps are market or idiosyncratic jumps. 

The DLE is fundamentally different than its continuous counterpart, which was studied by \cite{kx16} and \cite{aflwy16}, for example. First, in line with the model of \cite{pv12}, a negative DLE across stocks exists for market downward price jumps but DLE estimates are consistently positive for market upward price jumps. Second, financial leverage, measured by a firm's debt-to-equity ratio, only explains price-volatility cojump correlations for negative market jumps. It is the level of volatility which consistently explains the cross-section of DLE estimates.  

Our findings have implications for the parametric modeling of asset prices. Our empirical results cast doubt on the unconditional bivariate normality assumption of \cite{br16}, which implies tail independence and a generally linear relation around the center of the price-volatility cojump distribution. On the contrary, our results indicate that price-volatility cojumps around the center of the joint jump distribution\textemdash i.e.\ smaller jumps\textemdash are usually only weakly related, while jumps of the upper and the lower quantiles exhibit a strong and significant DLE. The linear dependence, which was introduced by \cite{dps00}, allows for tail dependence but imposes one linear relation for both upward and downward price jumps. This appears to be at odds with the data. A specification that combines the uncorrelatedness assumption of \cite{bcj07} and a price jump sign dependence, as modeled by \cite{mfm16}, appears as a candidate to adequately capture jump sizes of contemporaneous price and volatility cojumps. Working out the pricing implications of such a parametric model might be a path for future research.  

Finally, one would like to explore the cross-sectional and time series dimension of the estimated DLE in more detail. Since we discovered a significant link with the level of integrated volatility, it is natural to ask if an asset pricing framework, such as that in \cite{chw15}, prices discontinuous leverage.

\begin{appendix}
\renewcommand*{\thesection}{\appendixname}

\section{Proofs}
Standard localization techniques allow us to assume that there exists a constant $\Lambda$, such that
 \begin{align*}\max{\{|b_s(\omega)|,|\sigma_s(\omega)|,|X_s(\omega)|,|\delta_{\omega}(s,x)|/\gamma(x)\}}\le \Lambda\,,\end{align*}
for all $(\omega,s,x)\in(\Omega,\mathbb{R}_+,\mathbb{R})$; i.e., characteristics are uniformly bounded. We refer to \cite{jp12}, Section
4.4.1, for a proof.
\renewcommand*{\appendixname}{A}
\subsection{Proof of Proposition 3.1}
We decompose the observations $Y_{t_i^n}$ into signal $X_{t_i^n}$ and noise $\epsilon_{t_i^n}$. In order to analyze the discretization variance from the signal terms, an illustration of the pre-processed price estimates \eqref{locest} as a function in the efficient log-returns $\Delta_i^n X$ is helpful. Reordering addends, similar as in the proofs of \cite{z06}, we obtain the identity
\begin{align}\hspace*{-.35cm}\notag M_n^{-1}\Big(\sum_{i=l}^{l+M_n-1}Y_{t_i^n}-\sum_{i=l-M_n}^{l-1}Y_{t_i^n}\Big)&=M_n^{-1}\sum_{i=l}^{l+M_n-1}\big(Y_{t_i^n}-Y_{t_{i-M_n}^n}\big)\\
&\label{transformpa} =\sum_{k=1}^{M_n-1}\Delta_{l+k}^n Y\,\frac{M_n-k}{M_n}+\sum_{k=0}^{M_n-1}\Delta_{l-k}^nY\,\frac{M_n-k}{M_n}\,.\end{align}
The expectation and variance of noise terms are readily derived using the left-hand side of \eqref{transformpa} and the fact that $\epsilon_{t_i^n}$ is i.i.d.\,with mean zero and variance $\eta^2$. For the signal part, we exploit the above identity and consider the right-hand side of \eqref{transformpa}. Considering the drift part in the pre-processed price estimates \eqref{locest}, we can bound the right-hand side above by
\[\Big|\sum_{k=1}^{M_n-1}\Delta_{l+k}^n b\,\frac{M_n-k}{M_n}+\sum_{k=0}^{M_n-1}\Delta_{l-k}^nb\,\frac{M_n-k}{M_n}\Big|\le K\, M_n\,n^{-1}=\KLEINO(n^{-1/4})\,,\]
$\P$-almost surely, with a constant $K$, using the fact that 
\[\sum_{k=1}^{M_n-1} (1-k/M_n)+\sum_{k=0}^{M_n-1}(1-k/M_n)=M_n\,.\]
We decompose the signal process, $X_t=\int_0^t b_s\,ds+C_t+J_t$, into its jump component, $(J_t)_{t\ge 0}$, and the continuous It\^{o} semimartingale, $(C_t)_{t\ge 0}$.
Under Assumption \ref{jumps} and for $r\ge 1$, we can use the estimate
 \begin{align}\label{jumpestimate}\forall s,t\ge 0:~\E\left[|J_{t}-J_s|^p\big|\mathcal{F}_s\right]&\le K_p\,\E\Big[\Big(\int_s^{t}\int_{\R}(\gamma^r(x)\wedge 1)\lambda(dx)ds\Big)^{1/r}\Big]\\ \notag &  \le K_p |t-s|^{(1/r)}\, \end{align}
to find that the jump terms in the right-hand side of \eqref{transformpa} satisfy
\begin{align*}\E\Big[\Big|\sum_{k=1}^{M_n-1}\Delta_{l+k}^n J\,\frac{M_n-k}{M_n}+\sum_{k=1}^{M_n-1}\Delta_{l-k}^nJ\,\frac{M_n-k}{M_n}\Big|\Big]
=\mathcal{O}\big(M_n{n}^{-1/r}\big)\, ,\end{align*}
with some $r<4/3$, where we omit $\Delta_l^n J$ for $l=\lfloor \tau n\rfloor+1$. Thus, the terms multiplied with $n^{1/4}$ tend to zero in probability by Markov's inequality. Because the expectations of all increments $\Delta_i^n C$ vanish and $\E[\Delta_l^n (C+J+\epsilon)|\Delta X_{\tau}]=\Delta X_{\tau}$ for $l=\lfloor \tau n\rfloor+1$, we conclude that
\[\E\big[T^{LM}(\tau;\Delta_1^n Y,\ldots,\Delta_n^n Y)\big)|\Delta X_{\tau}\big]=\Delta X_{\tau}+\KLEINO_{\bar \P}\big(n^{-1/4}\big)\,.\]
In the case that $t_i^n=i/n$, It\^{o} isometry and the smoothness of the volatility granted by \eqref{JM2} and \eqref{JM} imply that for $l=\lfloor \tau n\rfloor+1$, 
\begin{align*}\E[(\Delta_{l+k}^n C)^2|\mathcal{F}_{\tau}]&=\E \Big[\int_{(l+k-1)/n}^{(l+k)/n}\sigma_s^2\,ds\big|\mathcal{F}_{\tau}\Big]+\mathcal{O}_{\P}(n^{-2})\\
&=\frac{\sigma_{\tau}^2}{n}+\mathcal{O}_{\P}\big(n^{-1} \sqrt{M_n/n}\big)\,,\end{align*} 
for all $k=1,\ldots,M_n-1$. Analogously, we obtain that
\begin{align*}\E[(\Delta_{l-k}^n C)^2|\mathcal{F}_{\tau-M_n/n}]&=\frac{\sigma_{\tau-}^2}{n}+\mathcal{O}_{\P}\big(n^{-1} \sqrt{M_n/n}\big)\,,\end{align*} 
for all $k=1,\ldots,M_n-1$.
Use of the identities
\begin{align*}
\sum_{k=0}^{M_n-1} (1-k/M_n)^2&=\frac13 M_n+\frac12+\frac16 M_n^{-1}\,,\\
\sum_{k=1}^{M_n-1} (1-k/M_n)^2&=\frac13 M_n-\frac12+\frac16 M_n^{-1}\,,
\end{align*}
and the independence of the noise and signal terms yield the asymptotic variance,
\[\var\big(\sqrt{M_n}\,T^{LM}(\tau;\Delta_1^n Y,\ldots,\Delta_n^n Y)\big)\rightarrow \frac13 (\sigma_{\tau}^2+\sigma_{\tau-}^2)\,c^2+2\eta^2\,\]
of the rescaled statistic. The form of the variance in \eqref{stablecltlm} follows from the above. 
Using that
\begin{align*}\E[(\Delta_{l+k}^n C)^4|\mathcal{F}_{\tau}]&=\frac{3\,\sigma_{\tau}^4}{n^2}+\KLEINO_{\P}\big(n^{-2}\big)~,~\E[(\Delta_{l-k}^n C)^4|\mathcal{F}_{\tau-M_n/n}]=\frac{3\,\sigma_{\tau-}^4}{n^2}+\KLEINO_{\P}\big(n^{-2}\big)\, \end{align*} 
and with the assumed existence of $\E[\epsilon_t^4]$, the Lyapunov criterion with fourth moments obtained from \eqref{transformpa} yields, together with the above considerations, the central limit theorem \eqref{stablecltlm}.

Next, we prove that the convergence is stable in law. The latter is equivalent to the joint weak convergence of $\alpha_n=\sqrt{M_n}\big(T^{LM}(\tau;\Delta_1^n Y,\ldots,\Delta_n^n Y)-\Delta X_{\tau}\big)$ with any $\mathcal{G}$-measurable bounded random variable $Z$:
\begin{align}\label{stable}\E\left[Z g(\alpha_n)\right]\rightarrow \E\left[Z g(\alpha)\right]=\E[Z]\E\left[g(\alpha)\right]\, \end{align}
for any continuous bounded function, $g$, and
\begin{align}\alpha=\big(1/3(\sigma_{\tau}^2+\sigma_{\tau-}^2)c^2+2\eta^2\big)^{1/2} U\,,\end{align}
with $U$ a standard normally-distributed, random variable that is independent of $\mathcal{G}$. 
In order to verify \eqref{stable}, consider the sequence $A_n=[(\tau-M_n/n)\vee 0, (\tau+M_n/n)\wedge 1]$.
Each $\alpha_n$ is measurable with respect to the $\sigma$-field $\mathcal{G}_1$. The sequence of decompositions 
\begin{align*}\tilde C(n)_t=\int_0^t\1_{A_n}(s)\sigma_{s}\,dW_s\,,\,\bar C(n)_t=C_t-\tilde C(n)_t\,,\end{align*}
\begin{align*}\tilde \epsilon(n)_t=\1_{A_n}(t)\epsilon_t\,,\,\bar \epsilon(n)_t=\epsilon_t-\tilde \epsilon(n)_t\,,\end{align*}
of $(C_t)_{t\ge 0}$ and $(\epsilon_t)_{t\ge 0}$ are well-defined. If $\mathcal{H}_n$ denotes the $\sigma$-field generated by $\bar C(n)_t$, $\bar \epsilon(n)_t$ and $\mathcal{F}_0$, then $\big(\mathcal{H}_n\big)_n$ is an isotonic sequence with $\bigvee_n \mathcal{H}_n=\mathcal{G}_1$. Since $\E[Z|\mathcal{H}_n]\rightarrow Z$ in $L^1(\P)$, it thus suffices that
\begin{align}\E[Z g (\alpha_n)]\rightarrow \E[Z\,g(\alpha)]= \E[Z]\E[g(\alpha)]\,,\end{align}
for $Z$ being $\mathcal{H}_q$ measurable for some $q$. Note that we can approximate the volatility to be constant over local intervals $[\tau-M_n/n,\tau)$ and $[\tau,\tau+M_n/n]$. Then, for all $n\ge q$, conditional on $\mathcal{H}_q$, $\alpha_n$ has a law independent of $\bar C(n)_t$ and $\bar \epsilon(n)_t$, such that the ordinary central limit theorem implies the claimed convergence.

\subsection{Proof of Proposition 3.2}
A neat decomposition of the spectral statistics into observation errors and returns of the efficient price is obtained with \emph{summation by parts}
\begin{align}\label{sbp}S_j(\tau)=\Big(\sum_{i=1}^n\Delta_i^n X\Phi_{j,\tau}((t_{i-1}^n+t_i^n)/2)-\sum_{i=1}^{n-1}\epsilon_{t_i^n}\Phi_{j,\tau}'(t_i^n)\,\frac{t_{i+1}^n-t_{i-1}^n}{2}\Big)(1+\KLEINO_{\bar\P}(1))\,,\end{align}
where the asymptotically negligible remainder comes from approximating\\ \(\Phi_{j,\tau}((t_{i+1}^n+t_i^n)/2)-\Phi_{j,\tau}((t_{i-1}^n+t_i^n)/2)\) with the derivative and end-effects. The system of derivatives $(\Phi_{j,\tau}')_{j\ge 1}$ is again orthogonal such that covariances between different spectral frequencies vanish.

First, we prove that the drift is asymptotically negligible under Assumption \ref{sigma}. Because $\int_0^1\Phi_{j,\tau}(t)\,dt=2\sqrt{2h_n}/(\pi j)$ and $\int_0^1\big|\Phi_{j,\tau}(t)\big|\,dt=2\sqrt{2h_n}/\pi$, we get with generic constant $K$ that $\P$-almost surely
\begin{align*}\Big|\sum_{i=1}^n\Delta_i^n b \,\Phi_{j,\tau}((t_{i-1}^n+t_i^n)/2)\Big|&\le K \sum_{i=1}^n(t_i^n-t_{i-1}^n)\big|\Phi_{j,\tau}((t_{i-1}^n+t_i^n)/2)\big|\le K\,\frac{\sqrt{h_n}}{\pi}\,,\end{align*} 
and thus
\begin{align*}\Big|\sum_{j=1}^{J_n}(-1)^{j+1}a_{2j-1}\sum_{i=1}^n\Delta_i^n b \,\Phi_{j,\tau}((t_{i-1}^n+t_i^n)/2)\Big|&\le K\sum_{j=1}^{J_n}(1+j^2h_n^{-2}/n)^{-1}\sqrt{h_n}\\
&=K\sum_{j=1}^{J_n}\Big(1+\frac{j^2}{\kappa^2\log^2(n)}\Big)^{-1}\sqrt{h_n}\\
&\le K\Big(\sum_{j=1}^{\log(n)}\sqrt{h_n}+\sum_{j=1}^{J_n}j^{-2}\sqrt{h_n}\log^2{(n)}\Big)\\
&\le K\log^2(n)\sqrt{h_n}\,.\end{align*}
This yields that $\P$-almost surely 
\[n^{1/4}\sqrt{\frac{h_n}{2}}\Big|\sum_{j=1}^{J_n}(-1)^{j+1}a_{2j-1}\sum_{i=1}^n\Delta_i^n b \,\Phi_{j,\tau}((t_{i-1}^n+t_i^n)/2)\Big|\rightarrow 0 \,,\]
which ensures that we can neglect the drift in the asymptotic analysis of \eqref{bw}.

Next, we analyze the variance of \eqref{bw} with oracle optimal weights \eqref{weights}. A locally constant approximation of $\sigma_s,s\in[\tau-h_n/2,\tau)$ and $\sigma_s,s\in[\tau,\tau+h_n/2]$ is asymptotically negligible under Assumption \ref{sigma}. Based on \eqref{sbp}, using the fact that
\begin{align*}\int_{\tau-h_n/2}^{\tau}\Phi_{j,\tau}^2(t)\,dt=\int_{\tau}^{\tau+h_n/2}\Phi_{j,\tau}^2(t)\,dt=1/2\, ,\end{align*}
yields the following variances of spectral statistics:
\begin{align*}\var\big(S_j(\tau)\big)=\frac{1}{2}(\sigma_{\tau}^2+\sigma_{\tau_-}^2)+\frac{\pi^2 j^2}{h_n^2}\frac{\eta^2}{n}\,.\end{align*}
We thus obtain the conditional variance,
\begin{align*}&\var\Big(n^{1/4}\,\mathcal{T}(\tau;\Delta_1^n Y,\ldots,\Delta_n^n Y)\big|\mathcal{F}_{\tau}\Big) \\ &=n^{1/2}\Big(\sum_{j=1}^{J_n}\big(\frac12(\sigma_{\tau}^2+\sigma_{\tau-}^2)+\pi^2(2j-1)^2h_n^{-2}n^{-1}\eta^2\big)^{-1}\Big)^{-1}h_n/2+\KLEINO_{\bar\P}(1)\\
&=\frac12\bigg(\frac{\sum_{j=1}^{J_n}(\frac12(\sigma_{\tau}^2+\sigma_{\tau-}^2)+\pi^2(2j-1)^2h_n^{-2}n^{-1}\eta^2)^{-1}}{\log(n)}\bigg)^{-1}+\KLEINO_{\bar\P}(1)\\
&=\frac12\,\Big(\int_0^{\infty}\frac{1}{\frac12(\sigma_{\tau}^2+\sigma_{\tau-}^2)+\pi^2(2z)^2\eta^2}\,dz\Big)^{-1}(1+\KLEINO(1))+\KLEINO_{\bar\P}(1)\\
&=2\Big(\frac{\sigma_{\tau}^2+\sigma_{\tau-}^2}{2}\Big)^{1/2}\eta+\KLEINO_{\bar\P}(1)\,.
\end{align*}
With $\delta_n\le n^{-1}$, $l=\lfloor \tau n\rfloor +1$, we have that
\begin{align*}&\E\big[\mathcal{T}(\tau;\Delta_1^n Y,\ldots,\Delta_n^n Y)|\Delta X_{\tau}\big]\\
&\qquad =\sqrt{\frac{h_n}{2}}\sum_{j=1}^{J_n}a_{2j-1}(-1)^{j+1}\Phi_{2j-1,\tau}(\tau+\delta_n)\E[\Delta_l^n Y|\Delta X_{\tau}]+\KLEINO_{\bar\P}(n^{-1/4})\\
&\qquad =\sqrt{\frac{h_n}{2}}\sum_{j=1}^{J_n}a_{2j-1}(-1)^{j+1}\Phi_{2j-1,\tau}(\tau+\delta_n)\Delta X_{\tau}+\KLEINO_{\bar\P}(n^{-1/4})\\
&\qquad =(1+\KLEINO(\delta_n))\Delta X_{\tau}+\KLEINO_{\bar\P}(n^{-1/4})\,.\end{align*}
Considering further jumps on the estimation window, utilizing \eqref{jumpestimate} yields 
\begin{align*}\E\bigg[\Big|\sum_{j=1}^{J_n}(-1)^{j+1}a_{2j-1}\sum_{ i\ne l}\Delta_i^n J \,\Phi_{j,\tau}((t_{i-1}^n+t_i^n)/2)\Big|\bigg]\le K\log^2(n)\sqrt{h_n}\sup_i|t_i^n-t_{i-1}^n|^{1/r-1}\end{align*}
by the triangle inequality, decomposing $|t_i^n-t_{i-1}^n|^{1/r}=(t_i^n-t_{i-1}^n)|t_i^n-t_{i-1}^n|^{1/r-1}$ and using the same Riemann sum approximation as for the drift terms above. As for the Lee-Mykland statistic, $r<4/3$ ensures asymptotic negligibility of further jumps on $[\tau-h_n/2,\tau+h_n/2]$.

Since we assume $\E[\epsilon_t^4]<\infty$, we can establish a Lyapunov condition with fourth moments. Integral approximations with $\int_0^{1}\Phi_{j,\tau}^4(t)\,dt$ and $\int_{0}^1(\Phi_{j,\tau}')^4(t)\,dt$ yield, with generic constant $C$, for all $j$,
\[n\frac{h_n^2}{4}\sum_{i=1}^n\E[(\Delta_i^n X)^4]\Phi_{j,\tau}^4((t_{i-1}^n+t_i^n)/2)\le C n h_n^2\,n^{-1}\,\frac32 h_n^{-1}=\mathcal{O}(h_n)\, \]
as well as,
\begin{align*}n\frac{h_n^2}{4}\sum_{i=1}^n\E[(\epsilon_{t_i^n})^4](\Phi_{j,\tau}')^4(t_i^n)(t_{i+1}^n-t_{i-1}^n)^4/16&\le C n h_n^2\,h_n^{-5}n^{-3}\log^5(n)\\
&\le C n^{-1/2}\log^2(n)\,.\end{align*}
Considering signal and noise terms separately, we derive for the signal terms with Jensen's inequality for weighted sums,
\begin{align*}&n\frac{h_n^2}{4}\sum_{i=1}^n\E\Big[\Big(\Delta_i^n X\sum_{j=1}^{J_n}(-1)^{j+1}a_{2j-1}\Phi_{2j-1,\tau}((t_{i-1}^n+t_i^n)/2)\Big)^4\Big]\\
&\le n\frac{h_n^2}{4}\sum_{i=1}^n\E[(\Delta_i^n X)^4]\sum_{j=1}^{J_n}a_{2j-1}\Phi_{2j-1,\tau}^4((t_{i-1}^n+t_i^n)/2)=\mathcal{O}\Big(\sum_{j=1}^{J_n} a_{2j-1} \,h_n\Big)=\mathcal{O}(h_n)\,.
\end{align*} 
An analogous bound by Jensen's inequality for the noise terms implies the Lyapunov condition.

Stability of weak convergence is proved along the same lines as for Proposition \ref{proplm} and we omit the proof. It remains to show that
\begin{align}\label{ada1}\E\Big[\Big|\sum_{j=1}^{J_n}(-1)^{j+1}\big(\hat a_{2j-1}-a_{2j-1}\big)S_{2j-1}(\tau)\,\sqrt{h_n/2}\Big|\Big]=\KLEINO_{\bar\P}\big(n^{-1/4}\big)\, ,\end{align}
where $\hat a_{2j-1}$ denote the estimated oracle weights, to prove the claimed result. Using the triangle and H\"older's inequalities, we can bound the right-hand side of \eqref{ada1} by
\begin{align*}&\sum_{j=1}^{J_n}\E\Big[\big|\hat a_{2j-1}-a_{2j-1}\big|\,\big| S_{2j-1}(\tau)\big|\Big]\sqrt{\frac{h_n}{2}}\\
&\hspace*{3cm}\le \sum_{j=1}^{J_n}\Big(\E\Big[\big|\hat a_{2j-1}-a_{2j-1}\big|^2\Big]\,\E\Big[\big| S_{2j-1}(\tau)\big|^2\Big]\Big)^{1/2}\sqrt{\frac{h_n}{2}}\,.\end{align*}
In order to analyze the magnitude of the error of pre-estimating the weights, $\big|\hat a_{2j-1}-a_{2j-1}\big|^2$, we interpret \eqref{weights} as a function of the variables $\sigma_{\tau}^2$, $\sigma_{\tau-}^2$ and $\eta^2$. Differential calculus and the delta method yield the upper bound,
\begin{align*}&\sum_{j=1}^{J_n}K\Big(a_{2j-1}^2\big(\delta_n(\sigma_{\tau}^2)+\delta_n(\eta^2)\big)^2\,\var(S_{2j-1}(\tau))\Big)^{1/2}\sqrt{h_n/2}\\
&\hspace*{-.5cm} \le \sum_{j=1}^{J_n}K\,\delta_n(\sigma_{\tau}^2)\frac{\big(\var(S_{2j-1}(\tau))\big)^{-1/2}}{\sum_{u=1}^{J_n}\big(\var(S_{2u-1}(\tau))\big)^{-1}}\sqrt{h_n/2}\\
&\hspace*{-.5cm} \le \sum_{j=1}^{J_n}K \Big(1+\frac{j^2}{\log^2(n)}\Big)^{-1/2}\delta_n(\sigma_{\tau}^2)\sqrt{h_n/2}=\mathcal{O}\Big(\log^3(n)\delta_n(\sigma_{\tau}^2)\sqrt{h_n}\Big)=\KLEINO\big(n^{-1/4}\big),\end{align*}
for the right-hand side of \eqref{ada1} with generic constant $K$ and bounds $\delta_n(\sigma_{\tau}^2)\le K n^{-1/8}$ and $\delta_n(\eta^2)\le K n^{-1/2}$ for the errors of pre-estimating $\sigma_{\tau}^2$, $\sigma_{\tau-}^2$ and $\eta^2$ with \eqref{pre1} and \eqref{pre2}, respectively. This ensures \eqref{ada1} and completes the proof of Proposition \ref{propbw}.
\subsection{Proof of Proposition 3.3}
The proof reduces to generalizing the analysis of the asymptotic variance and fourth moments for a Lyapunov condition.
Consider the noise term on the right-hand side of \eqref{sbp} under $R$-dependent noise and for $t_i^n=F^{-1}(i/n)$.
The expectation still vanishes and the variance becomes the following:
\begin{align*}&\E\Bigg[\Big(\sum_{i=1}^{n-1}\epsilon_{t_i^n}\Phi_{j,\tau}'(t_i^n)\tfrac{t_{i+1}^n-t_{i-1}^n}{2}\Big)^2\Bigg]\\
&=\E\Bigg[\sum_{i=1}^{n-1}\hspace*{-.05cm}\epsilon_{t_i^n}^2(\Phi_{j,\tau}'(t_i^n))^2\big(\tfrac{t_{i+1}^n-t_{i-1}^n}{2}\big)^2\hspace*{-.1cm}+\hspace*{-.05cm}2\hspace*{-.05cm}\sum_{i=1}^{n-1}\sum_{u=1}^{R\wedge (n-i)}\hspace*{-.15cm}\epsilon_{t_i^n}\epsilon_{t_{i+u}^n} \hspace*{-.05cm}\Phi_{j,\tau}'(t_i^n)\Phi_{j,\tau}'(t_{i+u}^n)\tfrac{t_{i+1}^n-t_{i-1}^n}{2}\tfrac{t_{i+u+1}^n-t_{i+u-1}^n}{2}\hspace*{-.05cm}\Bigg]\\
&= \E\Bigg[\sum_{i=1}^{n-1}(\Phi_{j,\tau}'(t_i^n))^2\tfrac{t_{i+1}^n-t_{i-1}^n}{2}(F^{-1})'(\tau)n^{-1}\Big(\epsilon_{t_i^n}^2+\sum_{u=1}^{R\wedge (n-i)}\epsilon_{t_i^n}\epsilon_{t_{i+1}^n}\Big)\Bigg](1+\KLEINO(1))\\
&=\eta^2_{\tau}(F^{-1})'(\tau)n^{-1}\int_0^1\Phi_{j,\tau}'(t)\,dt (1+\KLEINO(1))=\eta^2_{\tau}(F^{-1})'(\tau)n^{-1}\pi^2j^2h_n^{-2}(1+\KLEINO(1))\,.
\end{align*}
We used the smoothness of $(F^{-1})'$ and $\Phi_{j,\tau}'$ for approximations.
The same Riemann sum approximation as in the equidistant observations case applies for the signal term. Using a (double-Riemann sum) integral approximation as $J_n\rightarrow\infty$, analogously as in the proof of Proposition \ref{propbw}, yields the asymptotic variance in \eqref{stablecltbwgen}.
Introducing the shortcut, $\delta_{i,v}^R=\1_{\{|i-v|\le R\}}$, we obtain the following estimates for the fourth moments:
\begin{align*}&\E\Bigg[\Big(\sum_{i=1}^{n-1}\epsilon_{t_i^n}\Phi_{j,\tau}'(t_i^n)\tfrac{t_{i+1}^n-t_{i-1}^n}{2}\Big)^4\Bigg]\\
&=\E\Bigg[\hspace*{-.05cm}\sum_{i,v,u,r=1}^{n-1}\hspace*{-.05cm}\epsilon_{t_i^n}\epsilon_{t_v^n}\epsilon_{t_u^n}\epsilon_{t_r^n}\Phi_{j,\tau}'(t_i^n)\Phi_{j,\tau}'(t_v^n)\Phi_{j,\tau}'(t_u^n)\Phi_{j,\tau}'(t_r^n)\tfrac{t_{i+1}^{n}-t_{i-1}^{n}}{2}\tfrac{t_{v+1}^{n}-t_{v-1}^{n}}{2}\\
&\hspace*{8cm}\times\tfrac{t_{u+1}^{n}-t_{u-1}^{n}}{2}\tfrac{t_{r+1}^{n}-t_{r-1}^{n}}{2}\hspace*{-.05cm}\Bigg]\\
&=\sum_{i,v,u,r=1}^{n-1}\E\big[\epsilon_{t_i^n}\epsilon_{t_v^n}\epsilon_{t_u^n}\epsilon_{t_r^n} \big]\big(\delta_{i,v}^R\delta_{u,r}^R+\delta_{i,u}^R\delta_{v,r}^R+\delta_{i,r}^R\delta_{v,u}^R\big)\Phi_{j,\tau}'(t_i^n)\Phi_{j,\tau}'(t_v^n)\\
&\hspace*{3cm}\times \Phi_{j,\tau}'(t_u^n)\Phi_{j,\tau}'(t_r^n)\tfrac{t_{i+1}^{n}-t_{i-1}^{n}}{2}\tfrac{t_{v+1}^{n}-t_{v-1}^{n}}{2}\tfrac{t_{u+1}^{n}-t_{u-1}^{n}}{2}\tfrac{t_{r+1}^{n}-t_{r-1}^{n}}{2}\\
&=\big((F^{-1}(\tau))'\big)^23\,\eta^4_{\tau} \,n^{-2}-\mathcal{R}_n\,,\end{align*}
with a remainder, $\mathcal{R}_n$, that satisfies for some constant $C$ that
\begin{align*}&\mathcal{R}_n\le\sum_{i,v,u,r=1}^{n_p}C\,\Big(\delta_{i,v}^R\delta_{u,r}^R\big(\delta_{i,u}^R+\delta_{v,r}^R+\delta_{i,r}^R+\delta_{v,u}^R\big)+\delta_{i,u}^R\delta_{v,r}^R\big(\delta_{i,v}^R+\delta_{u,r}^R+\delta_{i,r}^R+\delta_{v,u}^R\big)\\
&\hspace*{6cm}+\delta_{i,r}^R\delta_{v,u}^R\big(\delta_{i,v}^R+\delta_{u,r}^R+\delta_{i,u}^R+\delta_{v,r}^R\big)\Big)n^{-4}\\
&=\mathcal{O}\big(nR^3n^{-4}\big)=\mathcal{O}\big(n^{-3}\big)=\KLEINO\big(n^{-2}\big)\,,\end{align*}
such that $\mathcal{R}_n$ is asymptotically negligible. Inserting the estimate, the Lyapunov condition is ensured in the generalized setting. Under $R$-dependence, the convergence of the generalized variance and the generalized Lyapunov criterion imply the central limit theorem \eqref{stablecltbwgen} and stability is proved analogously as above.
\subsection{Proof of Proposition 3.5}
Suppose that $\tau\in((k-1)h_n,kh_n)$ and we run the procedure from \eqref{argmax} to find a sub-interval that contains the jump.
The variances of the statistics $T^{LM}\big((k-1)h_n+(i-1/2)\frac{r_n+l_n}{n};\Delta_1^n Y,\ldots,\Delta_n^n Y\big),i=1,\ldots,R_n$, defined as in \eqref{lm} with $M_n$ replaced by $(r_n+l_n)/2$, are readily obtained from \eqref{transformpa} and given by
\begin{align*}&\var\Big(T^{LM}\big((k-1)h_n+(i-1/2)\frac{r_n+l_n}{n};\Delta_1^n Y,\ldots,\Delta_n^n Y\big)\Big)\\
&\hspace*{5cm}=\frac{4\eta^2_{(k-1)h_n}}{r_n+l_n}+\frac{l_n}{n}\frac{\sigma^2_{\tau-}}{3}+\frac{r_n}{n}\frac{\sigma^2_{\tau}}{3}+\KLEINO\Big(\frac{r_n+l_n}{n}\Big)\,.\end{align*}
In particular, for $r_n+l_n=\KLEINO(\sqrt{n})$, remainders in the proof of Proposition 3.1 become even smaller and the noise term prevails in the variance, such that, for all $i$, 
\[\sqrt{r_n+l_n}\,T^{LM}\big((k-1)h_n+(i-1/2)\tfrac{r_n+l_n}{n};\Delta_1^n (C+\epsilon),\ldots,\Delta_n^n (C+\epsilon)\big)\stackrel{(st)}{\longrightarrow}MN\big(0,4\eta^2_{(k-1)h_n}\big),\]
with $C_t$ the continuous semimartingale part of $X_t$. Since, under Assumption \ref{eta}, covariances between 
\[\Big(\frac{2}{r_n+l_n}\hspace*{-.1cm}\sum_{j=T_i}^{T_{i}+(r_n+l_n)/2-1}\hspace*{-.1cm}\big(\epsilon_{t_j^n}-\epsilon_{t^n_{j-(r_n+l_n)/2}}\big)\Big)_{i=1,\ldots,R_n}~,~T_i=\lfloor (k-1)h_n\,n\rfloor+(i-1/2)(r_n+l_n)+1,\]
are negligible, we deduce joint weak convergence to i.i.d.\;Gaussian limit random variables. Similarly as in \cite{lm12}, using basic extreme value theory, we derive that
\begin{align*}&B_n^{-1}\hspace*{-.05cm}\Big(\max_{i=1,\ldots,R_n}\hspace*{-.1cm}\sqrt{r_n+l_n}\,T^{LM}\big((k\hspace*{-.025cm} - \hspace*{-.05cm}1)h_n+(i\hspace*{-.025cm} - \hspace*{-.05cm}1/2)\tfrac{r_n+l_n}{n};\Delta_1^n (C+\epsilon),\ldots,\Delta_n^n (C+\epsilon)\big)\hspace*{-.05cm}-\hspace*{-.05cm}A_n\hspace*{-.05cm}\Big)\\
&\hspace*{12cm}\stackrel{(st)}{\longrightarrow}\xi\,,\end{align*}
with $\xi$ a standard Gumbel random variable  and
\[A_n=2\eta_{(k-1)h_n}\sqrt{2\log(R_n)}-\eta_{(k-1)h_n}\frac{\log(4\pi\log(R_n))}{\sqrt{2\log(R_n)}}~, ~ B_n^{-1}=\frac{\sqrt{\log(R_n)}}{\sqrt{2}\eta_{(k-1)h_n}}\,.\]
For $(r_n+l_n)\propto n^{\delta},\delta>0$, and if the jump is not located very close to the edges between the sub-intervals, the statistic on the sub-interval with the jump tends to infinity. That is, for
\[M_n=\max_{i=1,\ldots,R_n}\,T^{LM}\big((k-1)h_n+(i-1/2)\tfrac{r_n+l_n}{n};\Delta_1^n Y,\ldots,\Delta_n^n Y\big)~\, ,\]
where now $(\Delta_j^n Y)$ are inserted, not $(\Delta_j^n (C+\epsilon))$, we have that
\vspace*{-.4cm}

\[\sqrt{r_n+l_n}\,M_n-A_n\rightarrow\infty~\,.\]
As long as $M_n> (r_n+l_n)^{-1/2+\epsilon}$ for some $\epsilon>0$, this holds true. We need to carefully consider the potential bias issue discussed at the end of Section 3.1. The probability that $M_n\le (r_n+l_n)^{-1/2+\epsilon}$ translates to the probability that a jump is located in some small (in $n$ decreasing) vicinity of the block edges. Using that jump times are locally uniformly distributed, we obtain that 
\begin{align*}
\P\big(M_n\le (r_n+l_n)^{-\frac12+\epsilon}\big)&=\P\Big(\min_{i=1,\ldots,R_n}\Big|\tau-(k-1)h_n-(i-1)\tfrac{r_n+l_n}{n}\Big|\le \frac{(r_n+l_n)^{\frac12+\epsilon}}{n}\Big)\\
&=\P\big(U\in(0,2(r_n+l_n)^{-1/2+\epsilon})\big)=\mathcal{O}\big((r_n+l_n)^{-1/2+\epsilon}\big)\,,
\end{align*}
with $U$ a random variable uniformly distributed on $[0,1]$ and using the symmetry. Apparently, the probability converges to zero for $\epsilon$ sufficiently small. This implies that, for any such choice of $R_n$ and $r_n+l_n$, the procedure asymptotically almost surely detects the sub-interval which contains the jump.

Assigning jump times to a bin by thresholding induces a negligible error. This is proved analogously as in the proof of Proposition 3.6 in the next paragraph.

Cutting out noisy prices in the window $(t_{l-l_n}^n,t_{l+r_n}^n)$ around $\tau$, the adjusted statistics \eqref{bw} are asymptotically unbiased estimators of $\Delta X_{\tau}$. Considering their asymptotic properties, we can exploit most parts of the proof of Proposition \ref{propbw}. The only relevant difference is due to the increment over the cut-out window in the spectral statistics
\[\sqrt{\frac{h_n}{2}}\sum_{j=1}^{J_n}(-1)^{j+1}a_{2j-1}(Y_{t_{l+r_n}^n}-Y_{t_{l-l_n}^n})\Phi_{2j-1,\tau}(\tau)\,.\]
The increments $Y_{t_{l+r_n}^n}-Y_{t_{l-l_n}^n}$ take the role of $\Delta_l^n Y$, $l=\lfloor \tau n\rfloor+1$, where the window of statistics \eqref{bw} is centered. Using Jensen's inequality, we obtain that
\begin{align*}&\hspace*{-.25cm}\E\Big[\Big(\sqrt{\frac{h_n}{2}}\sum_{j=1}^{J_n}a_{2j-1}(-1)^{j+1}\Phi_{2j-1,\tau}(\tau)(X_{t_{l+r_n}^n}-X_{t_{l-l_n}^n})\Big)^2\Big]\\
&\le\frac{h_n}{2}\sum_{j=1}^{J_n} a_{2j-1} \Phi_{2j-1,\tau}^2(\tau)\E\big[(X_{t_{l+r_n}^n}-X_{t_{l-l_n}^n})^2\big]\\
&\le \max(\sigma_{\tau}^2,\sigma_{\tau-}^2)(t_{l+r_n}^n-t_{l-l_n}^n)=\KLEINO(n^{-1/2})\,.
\end{align*}
Since $\Phi_{2j-1,\tau}'(\tau)=0$, the summation by parts transformation \eqref{sbp} shows that the variance due to noise is not affected by the adjustment. Overall, we conclude that for the adjusted estimator,
\[|\widehat{\Delta X}_{\tau}-\mathcal{T}(\tau;\Delta_1^n Y,\ldots,\Delta_n^n Y)|=\KLEINO_{\bar\P}(n^{-1/4})\,.\]
We conclude the result with Proposition \ref{propbwgen}.
\subsection{Proof of Propositions 3.6 and 3.7}
Denoting  the finitely many stopping times with $|\Delta X_{\tau_k}|>a$, as $\tau_1,\ldots,\tau_N$, \eqref{dle} can be written 
\[[X,\sigma^2]^d_T(a)=\sum_{k=1}^N\Delta X_{\tau_k}\big(\sigma^2_{\tau_k}-\sigma^2_{\tau_k-}\big)\,.\]
The estimator \eqref{lev} then becomes 
\[\widehat{[X,\sigma^2]}^d_1(a)=\sum_{k=1}^{\hat N}\widehat{\Delta X}_{\hat \tau_k}\big(\hat\sigma^2_{\hat \tau_k}-\hat\sigma^2_{\hat \tau_k-}\big) \,.\]
The case without price jumps on the considered interval, $N=0$, is trivial. Consider the set
\begin{align*}\tilde \Omega_n=\left\{\omega\in\Omega|\tau_1>r_n^{-1}h_n,\tau_{N}<1-r_n^{-1}h_n,\tau_i-\tau_{i-1}>2r_n^{-1}h_n~,i=1,\ldots,N-1\right\}\\
\cup\; \left\{\omega\in\Omega|\tau_i= k\cdot h_n~,i=1,\ldots,N-1,k=0,\ldots,h_n^{-1}-1\right\}^{\complement}\,.\end{align*}
We can restrict to the subset $\tilde\Omega_n$, since $\P(\tilde\Omega_n)\rightarrow 1$ as $n\rightarrow \infty$. We infer the  jump times $\{\tau_i,i=1,\ldots,N\}$, or the respective bins on which jumps occur by thresholding. To show that this identification of jump times only induces an asymptotically negligible error, we prove that
\begin{align*}\Big|\sum_{k=2}^{h_n^{-1}-1}\widehat{\Delta X}_{\hat\tau_k}\big(\hat\sigma^2_{\hat\tau_k}-\hat\sigma^2_{\hat\tau_k-}\big)\,\1_{\{\Delta_k\widehat{[X,X]}>a^2 \vee u_n\}}-\sum_{k=1}^N\widehat{\Delta X}_{\tau_k}\big(\hat\sigma^2_{\tau_k}-\hat\sigma^2_{\tau_k-}\big)\Big|=\KLEINO_{\bar\P}\big(n^{-\beta/2}\big)\,.\end{align*}
This is ensured by Corollary \ref{testvolajump} and by Proposition \ref{argmaxprop} if 
\begin{align*}\sum_{k=2}^{h_n^{-1}-1}\big|\1_{\{\Delta_k\widehat{[X,X]}>a^2 \vee u_n\}}-\1_{\{\tau_i\in ((k-1)h_n,kh_n)\}}\big|=\KLEINO_{\bar\P}\big(n^{-1/8}\big)\,.\end{align*}
Denote $\mathcal{K}=\{1\le k\le h_n^{-1}|\tau_i\in ((k-1)h_n,kh_n)\}$ and $\mathcal{K}^{\complement}=\{2,\ldots, h_n^{-1}-1\}\setminus \mathcal{K}$. The last relation can be rewritten
\begin{align}\label{np}\sum_{k\in \mathcal{K}}\1_{\{\Delta_k\widehat{[X,X]}\le a^2 \vee u_n\}}+\sum_{k\in \mathcal{K}^{\complement}}\1_{\{\Delta_k\widehat{[X,X]}> a^2 \vee u_n\}}=\KLEINO_{\bar\P}\big(n^{-1/8}\big)\,.\end{align}
For each $k$ in the finite set $\mathcal{K}$, we prove that
\[\1_{\{\Delta_k\widehat{[X,X]}\le a^2 \vee u_n\}}=\KLEINO_{\bar\P}\big(n^{-1/8}\big)\,.\]
The restriction to $\tilde\Omega_n$ ensures that the considered jumps cannot occur on neighboring bins. Corollary 3.3 and its proof in \cite{bw15} establishes  that $\Delta_k\widehat{[X,X]}=(\Delta X_{\tau_i})^2+\chi_i$ with $\var(\chi_i)=\mathcal{O}(n^{-1/2})$. More precisely, as outlined in Section 3.1.3 of \cite{bw15}, for $\tau_i\in ((k-1)h_n,kh_n)$, we have that
\begin{align*}\E[h_n\,S_{jk}^2]&=2\sin^2(\pi j h_n^{-1}(\tau_i-(k-1)h_n))(\Delta X_{\tau_i})^2+\mathcal{O}(h_n)\,,\\
\E[h_n\,\max(\tilde S_{jk}^2,\tilde S_{j(k+1)}^2)]&=2\cos^2(\pi j h_n^{-1}(\tau_i-(k-1)h_n))(\Delta X_{\tau_i})^2+\mathcal{O}(h_n)\,.\end{align*}
The contribution with the cosine term is by $\tilde S_{jk}^2$ when $\tau_i\in ((k-1)h_n,(k-1/2)h_n)$ and by $ \tilde S_{j(k+1)}^2$ when $\tau_i\in ((k-1/2)h_n,kh_n)$. When $\tau_i=(k-1/2)h_n$, the cosine vanishes. Since the Lévy measure does not have an atom in $a$, it thus holds that, for some fixed $\epsilon>0$,
\[\Delta_k\widehat{[X,X]}=a^2+\epsilon+\chi_i\,.\]
Using Chebyshev's inequality, we derive that
\[\bar\P\big(\Delta_k\widehat{[X,X]}\le a^2 \vee u_n\big)\le \bar\P\big(|\chi_i|>\epsilon-u_n\big)=\mathcal{O}\big(n^{-1/2}\big)\,.\]
Considering indicator functions $\1_{A_n}$ with $p_n=\bar\P(A_n)\rightarrow 0$, using that $\E[\1_{A_n}]=p_n$ and $\var(\1_{A_n})\le p_n$, we obtain that
\begin{align}\label{np1}\sum_{k\in \mathcal{K}}\1_{\{\Delta_k\widehat{[X,X]}\le a^2 \vee u_n\}}=\mathcal{O}_{\bar\P}\big(n^{-1/4}\big)=\KLEINO_{\bar\P}\big(n^{-1/8}\big)\,.\end{align}
Due to the maximum operator in \eqref{zetatilde}, the term with the square cosine factor above feeds in two successive statistics. The cosine giving some factor bounded from above by one, we have for $\tau_i\in ((k-1)h_n,kh_n)$ that
\[\widetilde\zeta^{ad}_{k}>\max(\widetilde\zeta^{ad}_{k-1},\widetilde\zeta^{ad}_{k+1}) \,,\]
asymptotically almost surely. We conclude that the first sum in \eqref{np} is asymptotically negligible.

For $k\in \mathcal{K}^{\complement}$, neighboring a bin with $k\pm 1\in \mathcal{K}$, it holds asymptotically almost surely that $\widetilde\zeta^{ad}_{k}<\max(\widetilde\zeta^{ad}_{k-1},\widetilde\zeta^{ad}_{k+1})$, such that the indicator function sets it to zero. For all other $k\in \mathcal{K}^{\complement}$ we have that $\Delta_k\widehat{[X,X]}=h_n\tilde\zeta_k$ with $\tilde\zeta_k$ the local estimate for $\sigma^2_{(k-1)h_n}$ satisfying by Lemma 2 of \cite{bw16} the upper moment bound
\[\E\big[|\tilde\zeta_k|^{4+\delta}\big]=\mathcal{O}\big(\log(n)\big)\]
for $\delta $ from Assumption \ref{eta}. Markov's inequality yields for $k\in \mathcal{K}^{\complement}$
\[\bar\P\big(\Delta_k\widehat{[X,X]}> a^2 \vee u_n\big)=\mathcal{O}\Big(h_n^{4+\delta}\log(n)\big(a^2\vee u_n\big)^{(4+\delta)}\Big)\,.\]
Thereby we obtain that
\begin{align}\label{np2}\sum_{k\in \mathcal{K}^{\complement}}\1_{\{\Delta_k\widehat{[X,X]}> a^2 \vee u_n\}}=\mathcal{O}_{\bar\P}\Big(h_n^{-1}h_n^{2+\delta/2}\sqrt{\log(n)}\big(a^2\vee u_n\big)^{(2+\delta/2)}\Big)\,.\end{align}
If $a>0$, then this sum decays very fast as $n\rightarrow\infty$ and we clearly have $\KLEINO_{\bar\P}(n^{-1/8})$. If $a=0$, then the resulting order is $h_n^{1+\delta/2}\sqrt{\log(n)}h_n^{-\varpi(2+\delta/2)}$ and the term is $\KLEINO_{\bar\P}(n^{-1/8})$ when
\[\varpi<\frac{1+\delta/2-1/4}{2+\delta/2}\,.\]
A slightly stronger condition even implies the summability, 
\[\sum_{n\in\N}\bar\P\bigg(\sum_{k\in \mathcal{K}^{\complement}}\1_{\{\Delta_k\widehat{[X,X]}> a^2 \vee u_n\}}>0\bigg)<\infty\,,\]
and thus the almost sure convergence by Borel-Cantelli. With \eqref{np2}, we deduce \eqref{np} and are left to consider price and volatility jump estimates at times $\hat\tau_i,i=1,\ldots,\hat N$.
In both cases, $a>0$ or $a=0$ when $r=0$ in Assumption \ref{jumps}, $N<\infty$ holds almost surely. On $\tilde \Omega_n$, all involved local estimates for different price-jump times $\tau_i,i=1,\ldots,N$, are computed from disjoint datasets. The latter are not necessarily independent, but all covariations converge to zero asymptotically. For the single price-jump estimates, we have by Proposition 
\ref{argmaxprop} that 
\vspace*{-1cm}

\begin{subequations}
\begin{align*}
\widehat{\Delta X}_{\hat\tau_i}&=\Delta X_{\hat\tau_i}+\mathcal{O}_{\bar\P}\big(n^{-1/4}\big)\,.
\end{align*}
Based on Corollary \ref{testvolajump}, we obtain that
\begin{align*}
\Big(\big(\hat \sigma^2_{\hat\tau_i}-\hat \sigma^2_{\hat\tau_i -}\big)-\Delta\sigma^2_{\hat\tau_i}\Big)\stackrel{(st)}{\longrightarrow} \sqrt{8(\sigma_{\hat\tau_i}^3+\sigma_{\hat\tau_i-}^3)\eta_{\hat\tau_i}}Z_{i}\,,
\end{align*}
for all $i=1,\ldots,\hat N$, with $\beta$ from \eqref{beta} and $(Z_{i})$ i.i.d.\,standard normals. 
\end{subequations}
By the asymptotic negligibility of covariations, the vector 
\begin{align}\label{hvector}n^{\beta/2}\Big(\hat\sigma^2_{\hat\tau_1}-\hat\sigma^2_{\hat\tau_1-},\ldots,\hat\sigma^2_{\hat\tau_{\hat N}}-\hat\sigma^2_{\hat\tau_{{\hat N}}-}\Big)\stackrel{(st)}{\longrightarrow}\big(U_1,\ldots, U_{\hat N}\big)\,,\end{align}
converges stably in law, where $U_1,\ldots, U_{\hat N}$ are independent and \(U_i=\sqrt{8(\sigma_{\hat\tau_i}^3+\sigma_{\hat\tau_i-}^3)\eta_{\hat\tau_i}}Z_{i}\). 
Altogether, the asymptotic orders of different error terms and standard relations for weak and stochastic convergences imply $\eqref{cltdle}$. 
\subsection{Proof of Corollary 3.8}
According to the proof of Propositions \ref{propdle} and \ref{cordle}, the identification of bins with (large) jumps only induces an asymptotically negligible error. We are thus left to consider
\begin{align*}\frac{\sum_{k=1}^{\hat N}\widehat{\Delta X}_{\hat\tau_k}\big(\hat\sigma^2_{\hat\tau_k}-\hat\sigma^2_{\hat\tau_k-}\big)}{\Big(\sum_{k=1}^{\hat N}\big(\widehat{\Delta X}_{\hat\tau_k}\big)^2\sum_{k=1}^{\hat N}\big(\hat\sigma^2_{\hat\tau_k}-\hat\sigma^2_{\hat\tau_k-}\big)^2\Big)^{1/2}}-\frac{\sum_{k=1}^{\hat N}{\Delta X}_{\hat\tau_k}\Delta\sigma^2_{\hat \tau_k}}{\Big(\sum_{k=1}^{\hat N}\big({\Delta X}_{\hat\tau_k}\big)^2\sum_{k=1}^{\hat N}\big(\Delta\sigma^2_{\hat \tau_k}\big)^2\Big)^{1/2}}\,.\end{align*}
From \eqref{hvector}, we adopt that
\vspace*{-.3cm}

\[n^{\beta/2}\Big(\sum_{k=1}^{\hat N}\Big(\widehat{\Delta X}_{\hat\tau_k}\big(\hat\sigma^2_{\hat\tau_k}-\hat\sigma^2_{\hat\tau_k-}\big)-{\Delta X}_{\hat\tau_k}\Delta\sigma^2_{\hat \tau_k}\Big)\Big)\stackrel{(st)}{\longrightarrow}\sum_{k=1}^{\hat N}{\Delta X}_{\hat\tau_k} U_k\,.\]
\vspace*{-.3cm}

In the estimation of the discontinuous leverage, the estimation error of $(\Delta\sigma^2_{\hat \tau_k})$ dominates the smaller error of estimating $({\Delta X}_{\hat\tau_k})$. Analogously, for estimating \eqref{cor}
\vspace*{-.3cm}

\[\sum_{k=1}^{\hat N}\Big(\big(\widehat{\Delta X}_{\hat\tau_k}\big)^2-\big({\Delta X}_{\hat\tau_k}\big)^2\Big)=\mathcal{O}_{\bar \P}\big(n^{-1/4}\big)\,,\]
\vspace*{-.3cm}

readily obtained from Proposition \ref{propbwgen} and the delta method, induces an error that is negligible at first asymptotic order.
For the second variation, we deduce that
\vspace*{-.3cm}

\[n^{\beta/2}\sum_{k=1}^{\hat N}\Big(\big(\hat\sigma^2_{\hat\tau_k}-\hat\sigma^2_{\hat\tau_k-}\big)^2-\big(\Delta\sigma^2_{\hat \tau_k}\big)^2\Big)\stackrel{(st)}{\longrightarrow}\sum_{k=1}^{\hat N}2 \Delta\sigma^2_{\hat \tau_k}U_k\, ,\]
\vspace*{-.3cm}

from \eqref{hvector} and applying the binomial formula or the delta method for the square function. Another application of the delta method yields 
\begin{align*}&n^{\beta/2}\bigg(\frac{\widehat{[X,\sigma^2]}^d_1(a)}{\sqrt{\widehat{[X,X]}^d_1(a)\widehat{[\sigma^2,\sigma^2]}^d_1(a)}}-\frac{[X,\sigma^2]^d_1(a)}{\sqrt{[X,X]_1^d(a)[\sigma^2,\sigma^2]_1^d}(a)}\bigg)\\
&\quad =\sum_{k=1}^{\hat N_k} U_k\bigg(\frac{{\Delta X}_{\hat\tau_k}}{\big({[X,X]}^d_1(a){[\sigma^2,\sigma^2]}^d_1(a)\big)^{1/2}}-\frac{\Delta\sigma^2_{\hat \tau_k}[X,\sigma^2]^d_1(a)}{\big({[X,X]}^d_1(a)\big)^{1/2}\big({[\sigma^2,\sigma^2]}^d_1(a)\big)^{3/2}}\bigg)+\KLEINO_{\bar\P}(1)\,,\end{align*}
such that we obtain a stable central limit theorem with rate $n^{\beta/2}$ and asymptotic variance
\[\frac{8(\sigma_{\hat\tau_k}^3+\sigma_{\hat\tau_k-}^3)\eta_{\hat\tau_k}}{{[X,X]}^d_1(a){[\sigma^2,\sigma^2]}^d_1(a)}\Big({\Delta X}_{\hat\tau_k}-\frac{\Delta\sigma^2_{\hat \tau_k}[X,\sigma^2]^d_1(a)}{{[\sigma^2,\sigma^2]}^d_1(a)}\Big)^2\,.\]
This implies Corollary \ref{corcorrelation}.
\end{appendix}


\begin{thebibliography}{plainnat}
{\singlespacing 
\addcontentsline{toc}{section}{References}
\bibitem[A\"{i}t-Sahalia et al.(2017)]{aflwy16}A\"it-Sahalia,Y., Fan, J., Laeven, R.J.A., Wang, C.D.\ and Yang, X.\ (2017) \textit{Estimation of the continuous and discontinuous leverage effect}, forthcoming in Journal of the American Statistical Association.
\bibitem[A\"{i}t-Sahalia et al.(2013)]{afl13}A\"it-Sahalia,Y., Fan, J.\ and Li, Y.\ (2013) \textit{The leverage effect puzzle: Disentangling sources of bias at high frequency}, Journal of Financial Economics 109 (1), 224--249.
\bibitem[A\"{i}t-Sahalia and Jacod(2014)]{sahaliajacod}
{A\"{i}t-Sahalia}, Y. and {Jacod}, J.\ (2014) \textit{High-frequency financial econometrics}, Princeton, NJ: Princeton University Press.
\bibitem[A\"{i}t-Sahalia and Xiu(2017)]{ax16}
{A\"{i}t-Sahalia}, Y. and {Xiu}, D.\ (2017) \textit{A Hausman test for the presence of market microstructure noise in
high frequency data}, forthcoming in Journal of Econometrics.
\bibitem[Altmeyer and Bibinger(2015)]{ab15}Altmeyer, R. and Bibinger, M.\ (2015) \textit{Functional stable limit theorems for quasi-efficient spectral
covolatility estimators}, Stochastic Processes and their Applications 125 (12), 4556--4600.
\bibitem[Bandi and Ren\`o(2012)]{br12} Bandi, F.M.\ and Ren\`o, R. (2012) \textit{Time-varying leverage effects}, Journal of Econometrics 169 (1), 94--113.
\bibitem[Bandi and Ren\`o(2016)]{br16} Bandi, F.M.\ and Ren\`o, R. (2016) \textit{Price and volatility co-jumps}, Journal of Financial Economics 119 (1), 107--146.
\bibitem[Barndorff-Nielsen et al.(2009)]{bhls09}Barndorff-Nielsen, O. E., Hansen, P.R., Lunde, A.\ and Shephard, N.\ (2009) \textit{Realized kernels in practice: Trades and quotes}, Econometrics Journal 12 (3), C1--C32.
\bibitem[Bekaert and Wu(2000)]{bw00}Bekaert, G.\ and Wu, G.\ (2000) \textit{Asymmetric volatility and risk in equity markets}, Review of Financial Studies 13 (1), 1--42.
\bibitem[Benjamini and Hochberg(1995)]{bh95}Benjamini, Y.\ and Hochberg, Y.\ (1995) \textit{Controlling the false discovery rate: A practical and powerful approach to multiple testing}, Journal of the Royal Statistical Society Series B 57 (1), 289--300.
\bibitem[Bibinger et al.(2014)]{bhmr14}Bibinger, M., Hautsch, N., Malec, P.\  and Rei\ss, M.\  (2014) \textit{Estimating the quadratic covariation matrix from noisy observations: Local method of moments and efficiency}, Annals of Statistics 42 (4), 80--114.
\bibitem[Bibinger et al.(2017)]{bhmr17}Bibinger, M., Hautsch, N., Malec, P.\  and Rei\ss, M.\  (2017) \textit{Estimating the spot covariation of asset prices -- Statistical theory and empirical evidence}, forthcoming in Journal of Business and Economic Statistics.
\bibitem[Bibinger and Winkelmann(2015)]{bw15}Bibinger, M.\ and Winkelmann, L.\ (2015) \textit{Econometrics of cojumps in high-frequency data with noise}, Journal of Econometrics 184 (2), 361--378.
\bibitem[Bibinger and Winkelmann(2016)]{bw16}Bibinger, M.\ and Winkelmann, L.\ (2016) \textit{Common price and volatility jumps in noisy high-frequency data}, preprint, arxiv: 1407.4376.

\bibitem[Black(1976)]{b76}Black, F.\ (1976) \textit{Studies of stock price volatility changes}, Proceedings of the 1976 Meeting of the American Statistical Association, Business and Economic Statistics, 6, 177--181.

\bibitem[Bollerslev et al.(2006)]{blt06}Bollerslev, T., Litvinova, J.\ and Tauchen, G.\ (2006) \textit{Leverage and volatility feedback effects in high-frequency data}, Journal of Financial Econometrics 4 (3), 353--384.


\bibitem[Broadie et al.(2007)]{bcj07}Broadie, M., Chernov, M.\ and Johannes, M.\ (2007) \textit{Model specification and risk premia: Evidence from futures options}, Journal of Finance 62 (3), 1453--1490.
\bibitem[Chernov et al.(2003)]{cggt03}Chernov, M., Gallant, A.G., Ghysels, E.\ and Tauchen, G.\ (2003) \textit{Alternative models for stock price dynamics}, Journal of Econometrics 116 (1-2), 225--257.
\bibitem[Christensen et al.(2014)]{cop14}Christensen, K., Oomen, R.\ and Podolskij, M.\ (2014) \textit{Fact or friction: Jumps at ultra high frequency}, Journal of Financial Economics 114 (3), 576--599.
\bibitem[Cremers et al.(2015)]{chw15}Cremers, M., Halling, M.\ and Weinbaum, D.\ (2015) \textit{Aggregate jump and volatility risk in the cross-section of stock returns}, Journal of Finance 70 (2), 577--614.
\bibitem[Duffee(1995)]{d95}Duffee, G.R.\ (1995) \textit{Stock returns and volatility -- A firm-level analysis}, Journal of Financial Economics 37 (3), 399--420.
\bibitem[Duffie et al.(2000)]{dps00}Duffie, D., Pan, J.\ and Singleton, K.\ (2000) \textit{Transform analysis and asset pricing for affine-diffusions}, Econometrica 68 (6), 1343--1376.
\bibitem[Eraker(2004)]{e04}Eraker, B.\ (2004) \textit{Do stock prices and volatility jump? Reconciling evidence from spot and option prices}, Journal of Finance 59 (3), 1367--1403.
\bibitem[Eraker et al.(2003)]{ejp03}Eraker, B., Johannes, M.\ and Polson, N.\ (2003) \textit{The impact of jumps in volatility and returns}, Journal of Finance 58 (3), 1269--1300.
\bibitem[French et al.(1987)]{fss87}French, K.R., Schwert, W.\ and Stambaugh, R.\ (1987) \textit{Expected stock returns and volatility}, Journal of Financial Economics 19 (1), 3--29.
\bibitem[Hansen and Lunde(2006)]{hans06} Hansen, P.~R. and Lunde, A.\ (2006) \textit{Realized variance and market microstructure noise}, Journal of Business and Economic Statistics 24 (2), 127--161.
\bibitem[Jacod et al.(2017)]{jkm13}Jacod, J., Kl\"uppelberg, C.\ and M\"uller, G.\ (2017) \textit{Testing for non-correlation between price and volatility jumps}, Journal of Econometrics 197 (2), 284--297.
\bibitem[Jacod et al.(2010)]{jpv10}Jacod, J., Podolskij, M.\ and Vetter, M.\ (2010) \textit{Limit theorems for moving averages of
discretized processes plus noise}, Annals of Statistics 38 (3), 1478--1545.
\bibitem[Jacod and Protter(2012)]{jp12}Jacod, J.\ and Protter, P.\ (2012) \textit{Discretization of processes}, Springer.
\bibitem[Kalnina and Xiu(2017)]{kx16}Kalnina, I.\ and Xiu, D.\ (2017) \textit{Nonparametric estimation of the leverage effect: A trade-off between robustness and efficiency}, Journal of the American Statistical Association 112 (517), 384--396.
\bibitem[Koike(2017)]{Koike2017} Koike, Y.\ (2017) \textit{Time endogeneity and an optimal weight function in pre-averaging covariance estimation}, Statistical Inference for Stochastic Processes 20 (1), 15--56.
\bibitem[Lahaye et al.(2011)]{lln11}Lahaye, J., Laurent, S.\ and Neely, C.J.\ (2011) \textit{Jumps, cojumps and macro announcements}, Journal of Applied Econometrics 26 (6), 893--921.
\bibitem[Lee and Mykland(2008)]{lm8}Lee, S.\ and Mykland, P.A.\  (2008) \textit{Jumps in financial markets: A new nonparametric test and jump dynamics}, Review of Financial Studies 21 (6), 2535--2563.
\bibitem[Lee and Mykland(2012)]{lm12}Lee, S.\ and Mykland, P.A.\ (2012) \textit{Jumps in equilibrium prices and market microstructure noise}, Journal of Econometrics 168 (2), 396--406.
\bibitem[Li et al.(2017)]{ltt17}Li, J., Todorov, V.\ and Tauchen, G.\ (2017) \textit{Jump regressions}, Econometrica 85 (1), 173--195.
\bibitem[Maneesoonthorn et al.(2017)]{mfm16}Maneesoonthorn, W., Forbes, C.S.\ and Martin, G.M.\ (2017) \textit{Inference on self-exciting jumps in prices and volatility using high frequency measures}, Journal of Applied Econometrics 32 (3), 504--532.
\bibitem[Mykland and Zhang(2016)]{mykzhang16}Mykland, P.A. and Zhang, L.\ (2016) \textit{Between data cleaning and inference: Pre-averaging and robust estimators of the efficient price}, Journal of Econometrics 194 (2), 242--262.
\bibitem[P\'astor and Veronesi(2012, 2013)]{pv12}P\'astor, L.\ and Veronesi, P.\ (2012) \textit{Uncertainty about government policy and
stock prices}, Journal of Finance 67 (4), 1219--1264.
\bibitem[P\'astor and Veronesi(2013)]{pv13}P\'astor, L.\ and Veronesi, P.\ (2013) \textit{Political uncertainty and risk premia}, Journal of Financial Economics 110 (3), 520--545.
\bibitem[{Pelger(2017)}]{p17}Pelger, M.\ (2017) \textit{Understanding systematic risk: A high-frequency approach}, working paper, Department of Management Science \& Engineering, Stanford University, Stanford, CA 94305.
\bibitem[{Rei\ss(2011)}]{reiss}Rei\ss, M.\ (2011) \textit{Asymptotic equivalence for inference on the volatility from noisy observations}, Annals of Statistics 39 (2), 772--802.
\bibitem[Todorov and Tauchen(2011)]{tt11}Todorov, V.\ and Tauchen, G.\ (2011) \textit{Volatility jumps}, Journal of Business and Economic Statistics 29 (3), 356--371.
\bibitem[Vetter(2012)]{v12}Vetter, M.\ (2012) \textit{Estimation of correlation for continuous semimartingales}, Scandinavian Journal of Statistics 39 (4), 757--771.
\bibitem[Vetter(2014)]{v14}Vetter, M.\ (2014) \textit{Inference on the L\'{e}vy measure in case of noisy observations}, Statistics and Probability Letters 87, 125--133.
\bibitem[Wang and Mykland(2014)]{wm14}Wang, C.D.\ and Mykland, P.A.\ (2014) \textit{The estimation of the leverage effect with high-frequency data}, Journal of the American Statistical Association 109 (505), 197--215.
\bibitem[Winkelmann et al.(2016)]{wbl16}Winkelmann, L., Bibinger, M.\ and Linzert, T.\ (2016) \textit{ECB monetary policy surprises: Identification through cojumps in interest rates}, Journal of Applied Econometrics 31 (4), 613--629.
\bibitem[Yu(2012)]{y12} Yu, J.\ (2012) \textit{A semiparametric stochastic volatility model}, Journal of Econometrics 167 (2), 473--482.
\bibitem[Zhang(2006)]{z06} Zhang, L.\ (2006) \textit{Efficient estimation of stochastic volatility using noisy observations: A multi-scale approach}, Bernoulli 12 (6), 1019--1043.

}

\end{thebibliography}
\end{document}